\documentclass{amsart}
\usepackage{amssymb,amscd}
\usepackage{enumerate}
\begin{document}
\setcounter{tocdepth}{1}

%%%%%%%%%%%%%%%%%%%%%%%%%
% HAVE A NICE ``DATE:'' %
%%%%%%%%%%%%%%%%%%%%%%%%%
% The following command shows the labels.
%\let\labelc\label\renewcommand\label[1]{\mar{#1}\labelc{#1}}
%
%\setlength{\marginparwidth}{1.12in}
%\newcommand{\mar}[1]{{\marginpar{\textsf{#1}}}}
%\newcommand\datver[1]{\def\datverp%
%{\par\boxed{\boxed{\text{Version: #1; Run: \today}}}}}
%\datver{2.0; Revised:\ September 09, 2002 by V.}

%
% OPERATORS:
%
\newcommand\coker{\operatorname{coker}}
\newcommand\hotimes{\hat \otimes}
\newcommand\ind{\operatorname{ind}}
\newcommand\End{\operatorname{End}}
\newcommand\pa{\partial}
\newcommand\sign{\operatorname{sign}}
\newcommand\supp{\operatorname{supp}}
\newcommand\CI{\mathcal{C}^\infty}
\newcommand\CIc{\mathcal{C}^\infty_c}
\newcommand\CO{\mathcal{C}_0}
\newcommand\Hom{\operatorname{Hom}}
\newcommand\END{\operatorname{END}}
\newcommand{\dist}{\operatorname{dist}}
\newcommand\Dir{{\not \!\!D}} %%% old: {D\!\!\!\!\!\!/}
\newcommand{\Cliff}{{\rm Cliff}}
\newcommand\adj{\operatorname{ad}}
\newcommand\<{\langle}
\renewcommand\>{\rangle}
\newcommand{\ha}{{\textstyle{\frac12}}}

%
% LAYER POTENTIAL INTEGRALS
%

\newcommand\Sl{\mathcal{S}}
\newcommand\Dl{\mathcal{D}}

%
% PSEUDODIFFERENTIAL OPERATORS:
%
\newcommand{\PS}[1]{\Psi^{#1}(M)}
\newcommand{\IPS}[1]{\Psi_{\rm inv}^{#1}(M)}
\newcommand{\APS}[1]{\Psi_{\rm ai}^{#1}(M)}
\newcommand{\APSX}[1]{\Psi_{\rm ai}^{#1}(M \times X)}
\newcommand{\IPSN}[1]{\Psi_{\rm inv}^{#1}(\pa N)}
\newcommand{\APSN}[1]{\Psi_{\rm ai}^{#1}(\pa N)}
\newcommand{\IPSZ}[1]{\Psi_{\rm inv}^{#1}(Z)}
\newcommand{\APSZ}[1]{\Psi_{\rm ai}^{#1}(Z)}

%
% BOLD FACE SYMBOLS:
%
\newcommand\CC{\mathbb C}
\newcommand\NN{\mathbb N}
\newcommand\RR{\mathbb R}
\newcommand\ZZ{\mathbb Z}

%
% GENERAL LATEX MACROS:
%
\newcommand\ie{{i.\thinspace e.{}, }}
\newcommand\mfk{\mathfrak}

%
% THEOREM TYPE ENVIRONMENTS:
%
\newtheorem{theorem}{Theorem}[section]
\newtheorem{proposition}[theorem]{Proposition}
\newtheorem{corollary}[theorem]{Corollary}
\newtheorem{conjecture}[theorem]{Conjecture}
\newtheorem{lemma}[theorem]{Lemma}
\newtheorem{notation}[theorem]{Notations}
\theoremstyle{definition}
\newtheorem{definition}[theorem]{Definition}
\theoremstyle{remark}
\newtheorem{remark}[theorem]{Remark}
\newtheorem{example}[theorem]{Example}

%%%%%%%%%%%%%%%%%%%%%%%%%%%%%%%%%
%%                              %
%%  MACROS FOR THIS MANUSCRIPT  %
%%                              %
%%%%%%%%%%%%%%%%%%%%%%%%%%%%%%%%%
\newcommand\alge{\mfk A}
\newcommand\manif{\mfk M}

\author[M. Mitrea]{Marius Mitrea} \address{University of Missouri-Columbia,
        Department of Mathematics, Columbia, MO 65211}
        \email{marius@math.missouri.edu}

\author[V. Nistor]{Victor Nistor} \address{Pennsylvania State
       University, Math. Dept., University Park, PA 16802}
       \email{nistor@math.psu.edu}

\thanks{Mitrea was partially supported by NSF Grant DMS-0139801 and a
UMC Research Board Grant. Nistor was partially supported by NSF Grants
DMS-0209497 and DMS-0200808. Manuscripts available from {\bf http:
{\scriptsize//}www.math.missouri.edu{\scriptsize/}\,$\tilde{\,}$\,marius}
and {\bf
http:{\scriptsize//}www.math.psu.edu{\scriptsize/}nistor{\scriptsize/}}.
\\ \mbox{}\quad{2000 {\it Mathematics Subject Classification:} Primary
31C12, 58J05; Secondary 58J40, 35J05\qquad\qquad\qquad\qquad {\it Key
words and phrases}: layer potentials, manifolds with cylindrical ends,
Dirichlet problem}}

\date{\today}

%\dedicatory\datverp
\begin{abstract}\
We extend the method of layer potentials to manifolds with boundary
and cylindrical ends. To obtain this extension along the classical
lines, we have to deal with several technical difficulties due to the
non-compactness of the boundary, which prevents us from using the
standard characterization of Fredholm and compact
(pseudo-)differential operators between Sobolev spaces.  Our approach,
which involves the study of layer potentials depending on a parameter
on compact manifolds as an intermediate step, yields the invertibility
of the relevant boundary integral operators in the global, non-compact
setting, which is rather unexpected.  As an application, we prove the
well-posedness of the non-homogeneous Dirichlet problem on manifolds
with boundary and cylindrical ends.  We also prove the existence of
the Dirichlet-to-Neumann map, which we show to be a pseudodifferential
operator in the calculus of pseudodifferential operators that are
``almost translation invariant at infinity,'' a calculus that is
closely related to Melrose's b-calculus \cite{me81, meaps}, which we
study in this paper. The proof of the convergence of the layer
potentials and of the existence of the Dirichlet-to-Neumann map are
based on a good understanding of resolvents of elliptic operators that
are translation invariant at infinity.
\end{abstract}

\title[Manifolds with cylindrical ends]{Boundary value problems
and layer potentials on manifolds with cylindrical ends}

\maketitle \tableofcontents

\section*{Introduction}\
Boundary value problems, mostly on compact manifolds, have long
been studied because of their numerous applications to other areas
of Mathematics, Physics, and Engineering. Arguably, some of the
most important examples arise in connection with the Laplacian and
related operators.

A first, simple approach to boundary value problems for the
Laplace operator is via the Lax-Milgram theorem which amounts to
proving an energy estimate (coercivity) for the de Rham
differential of certain classes of scalar functions. Another
approach commonly used in the literature is via boundary layer
potential integrals. While less elementary, this has the advantage
that it provides more information about the spaces of Cauchy data,
and it allows one to express the solutions via explicit formulas.
The same approach may be used to study boundary problems on spaces
with weights, on which the Laplace operator may fail to be symmetric.

A second approach, based on the method of layer potentials, became 
widely used after the pioneering work of Hodge, de Rham, Kodaira, Spencer,
Duff, and Kohn, among others. See, for instance, \cite{DeR, Ho, Ko, KS}, 
or the discussion in the introduction of \cite{MMT} for further information 
and references. This method, combining ideas both from the
approach based on the Lax-Milgram theorem and the approach based on
the Boutet de Monvel calculus, has been successfully employed 
to solve boundary value problems on compact manifolds with smooth
boundary. 

More recently, the method of layer potentials has also lead to a
solution of the Dirichlet problem for the Laplace operator on compact
manifolds with Lipschitz boundaries in \cite{MTJFA}.  This, in turn,
builds on the earlier work from \cite{DK, FJR}, and \cite{Ve}, in the
constant coefficient, Euclidean context.

In view of several possible applications, some of which will be
discussed below, we would like to extend the method of layer
potentials to various classes of {\it non-compact} manifolds. There
are, however, several technical problems that we need to overcome
first for such an extension to be possible -- at least along the
classical lines.  The main contribution of this paper is to explain
how these difficulties are dealt with in the particular case of
manifolds with cylindrical ends, when a number of required results
from analysis take a simpler form. See also \cite{Grubb, ErkipSchrohe,
SchroheSG} for earlier results on boundary value problems on
non-compact manifolds.

The crucial step is to prove the invertibility of $-1/2 I + K$,
where $K$ is a suitable pseudodifferential operator, which is a
surprising, yet fortunate result. We hope that the approach that
we outline in this paper will serve as a paradigm for treating
more general elliptic PDE's on non-compact manifolds.

In fact, in  \cite{Schrohe2}, Elmar Schrohe has studied boundary
value problems for ``asymptotically Euclidean manifolds'' (this is
a class of non-compact manifolds generalizing the class of
manifolds that are Euclidean at infinity, a class introduced by
Chochet-Bruhat and Christodoulou, \cite{ChBrCh}). For this class
of manifolds, he has generalized the Boutet de Monvel's algebra
\cite{Boutet}. He has also pointed out the importance and
relevance of the spectral invariance of various algebras of
pseudodifferential operators. A main analytic difference between
his class of manifolds and ours is that while the ``Fredholm
relevant symbol'' (in the terminology of Cordes \cite{Cordes}) is
commutative for asymptotically Euclidean manifolds, this is no
longer true in the case we intend to study, i.e. that of manifolds
with cylindrical ends. Manifolds with cylindrical ends have also appeared 
in the study of boundary value problems on manifolds with conical points
\cite{Kondratiev, MazyaBook}. The results of this paper were used to 
prove the well-posedness of the Dirichlet problem in suitable Sobolev 
spaces with weights in \cite{BNZ}. This well-posedness result was then 
used in the same paper to obtain fast algorithms for solving the Dirichlet 
problem on polygonal domains in the plane.

In order to explain some of the technical difficulties encountered
in the setting of manifolds with cylindrical ends, we need to
introduce some notation. Let $N$ be a non-compact Riemannian
manifold with boundary $\pa N$ and $\Delta_N = d^*d$ be the
Laplace operator on $N$ action on scalar functions. A first set of
problems consists of defining an elementary solution
$E(\,\cdot\,,\,\cdot\,)$ for $\Delta_N$ on $N$ and proving that
the associated single and double layer potential integrals
converge --\,issues well-understood when $\pa N$ is compact. A
second set of problems has to do with the existence of the
non-tangential limits of the aforementioned layer potential
integral operators. Even if the non-tangential limits exist and
are given by pseudodifferential operators on $\pa N$, these
pseudodifferential operators are not expected to be properly
supported. Moreover, since $\pa N$ is non-compact, the standard
results on the boundedness and compactness of order zero
(respectively, negative order) pseudodifferential operators do not
(directly) apply. Finally, on non-compact manifolds one is lead to
consider various algebras of pseudodifferential operators with a
controlled behavior at infinity. These algebras may fail to be
``spectrally invariant,'' in the sense that the inverse of an
elliptic, $L^2$-invertible operator in this algebra may fail to be
again in this algebra.

In order to make the above technical problems more tractable, it is
natural to make certain additional assumptions on the non-compact
manifolds $N$ and $\pa N$ and, in this paper, we have restricted
ourselves to the class of manifolds with boundary and cylindrical
ends. For the sake of this introduction, let us briefly discuss about
pseudodifferential operators in this setting and then describe our
main results.

Let $M$ be a {\em boundaryless} manifold with cylindrical ends. Such
manifolds have a product structure at infinity in a strong sense (that
is, including also the metric --\,see Definition
\ref{def.m.b.c.e}). In this setting, we define two classes of
pseudodifferential operators: $\IPS{m}$ and $\APS{m}$, whose
distribution kernels form a class large enough to contain the
distribution kernels appearing in our paper as boundary layer
integrals.  See also \cite{Lewis, LewisParenti}, where some of these
issues were studied in the case of a polygon.

The first class of operators is the class of order $m$ classical
pseudodifferential operators that are ``translation invariant in a
neighborhood of infinity'' (Definition \ref{def.t.i}).  The space
$\APS{-\infty}$ consists of the closure of $\IPS{-\infty}$ with
respect to a suitable family of semi-norms, including for example the
norms of linear maps between the Sobolev spaces $H^{m}(M) \to
H^{m'}(M)$, $m,m' \in 2\ZZ$ (see Equations \eqref{eq.def.H2m} and
\eqref{eq.def.fam.s}; Sobolev spaces on non-integral orders can also
be defined, but they are not needed to construct our algebras). Then
\begin{equation}\label{eq.def.ai}
    \APS{m} := \IPS{m} + \APS{-\infty}.
\end{equation}
An operator $P \in \APS{m}$ is called {\em almost invariant in a
neighborhood of infinity}. For $P \in \APS{m}$, we can characterize
when it is Fredholm or compact (between suitable Sobolev spaces),
along the classical lines. See \cite{Kondratiev, LMN, LockhartOwen, MeMe,
MeN, Schrohe, ScSc, Schulze} and others.

We could have also allowed a power law behavior at infinity for our operators.
However, this is technical and would have greatly increased the size of the
paper, without really making our results more general. It would have also
shifted the focus of our paper, which is on boundary value problems and
not on constructing and studying algebras of pseudodifferential operators.

The reason for introducing the algebras $\APS{\infty}$ is that $T^{-1}
\in \APS{-m}$, for any elliptic operator $T \in \APSN{m}$, provided
that $m \ge 0$ and $T$ is elliptic and invertible on $L^2(M)$. (Recall
that $T$ is invertible as an unbounded operator if $T$ is injective
and $T^{-1}$ extends to a bounded operator.) This allows us to define
our integral kernels--and implicitly also the boundary layer
integrals--much as in \cite{MTJFA}, namely as follows. First, we embed
our manifold with boundary and cylindrical ends $N$ into a
boundaryless manifold with cylindrical ends $M$.  We then prove that
for suitable $V \ge 0$, $V \not = 0$, the operator $\Delta_M + V$ is
invertible by checking that it is Fredholm of index zero and
injective.

The single layer potential integral is defined then as
\begin{equation}
    \Sl(f) = (\Delta_M + V)^{-1} (f \otimes \delta_{\pa N}),
\end{equation}
where $f \in L^2(\pa N)$ and $\delta_{\pa N}$ the conditional measure
on $\pa N$ (so that $f \otimes \delta_{\pa N}$ defines the distribution
$\langle f \otimes \delta_{\pa N},\varphi\rangle = \int_{\pa N}f\varphi$,
where $\varphi$ is a test function on $N$).
Similarly, he double layer potential integral is
defined as
\begin{equation}
    \Dl(f) = (\Delta_M + V)^{-1} (f \otimes \delta'_{\pa N}),
\end{equation}
where $f \in L^2(\pa N)$, again, and $\delta'_{\pa N}$ the normal derivative
of the measure $\delta_{\pa N}$ in the sense of distributions (so that
$\langle f \otimes \delta'_{\pa N},\varphi\rangle=\int_{\pa N}f
\partial_\nu \varphi$, where $\nu$ is the unit normal to $N$).
Since we are dealing with non-compact manifolds ($M$ and $\pa N$), the
above integrals are defined by relying on mapping properties of the
operators in $\APS{m}$.

Next, we show that we can make sense of the restriction to $\pa N$ of
the kernel $E$ of $(\Delta_{M} + V)^{-1}$ and that the restricted kernel
gives rise to an operator
\begin{equation}\label{eq.def.S0}
    S : = [(\Delta_{M} + V)^{-1}]_{\pa N} \in\APSN{-1}.
\end{equation}

We can then relate the non-tangential limits of the single and double
layer potentials of some function $f$ using the operator $S$. This is
proved by writing $(\Delta_M + V)^{-1}$ as a sum of an operator $P \in
\IPSN{-1}$ and an operator $R \in \APS{-\infty}$. The existence and
properties of the integrals defined by $P$ follow as in the classical
case, because $P$ is properly supported (and hence all our relations
can be reduced to the analogous relations on a compact manifold). The
existence and properties of the integrals defined by $R$ follow from
the fact that $R$ is given by a uniformly smooth kernel, albeit not
properly supported.

We shall fix in what follows a vector field $\pa_\nu$ on $M$ that is
normal to $N$ at every point of $N$. Similarly, we define
\begin{equation}\label{eq.def.K0}
    K : = [(\Delta_{M} + V)^{-1} \pa_\nu^*]_{\pa N},
\end{equation}
by restricting the kernel of $(\Delta_{M} + V)^{-1} \pa_\nu^*$ to $\pa N$.
Let $f_{\pm}$ be the non-tangential pointwise limits  of some function
$f$ defined on $M \smallsetminus \pa N$, provided that they exist.

Some of the properties of the single and double layer potentials alluded to
above are summarized in the following theorem.

\begin{theorem}\label{theorem.map.sdlp2}\
Given $f \in L^2(\pa N)$, we have
\begin{equation}
\begin{gathered}
    \Sl(f)_+ = \Sl(f)_- = S f\,,\qquad \pa_\nu \Sl(f)_{\pm} =
    \Bigl(\pm \ha I + K^* \Bigr)f\, \quad \text{and}\\
    \Dl(f)_{\pm} = \Bigl(\mp \ha I + K\Bigr) f,
\end{gathered}
\end{equation}
where $K^*$ is the formal transpose of $K$.
\end{theorem}

These theorems are proved by reduction to the compact case
\cite{MTJFA} (using the decomposition $(\Delta_M + V)^{-1} = P + R$
explained above).  As in the classical case of a compact manifold with
smooth boundary, we obtain the following result.

\begin{theorem}\label{theorem.main}\
Let $N$ be a manifold with boundary and cylindrical ends. Then
$$
    H^{s}(N) \ni u \mapsto (\Delta_N u, u\vert_{\pa N})
    \in H^{s-2}(N) \oplus H^{s-1/2}(\pa N)
$$
is a continuous bijection, for any $s>1/2$.
\end{theorem}

A possible application of our results on boundary value problems
on manifolds with cylindrical ends is to Gauge theory, where
manifolds with cylindrical ends are often used. There are several
other potential applications of our results and methods to index
theory and to spectral theory on non-compact manifolds, not
necessarily with cylindrical ends.  See \cite{meaps, Muller1,
Muller2}. When extended to Dirac operators, our results, we hope,
will be useful to study Hamiltonians whose potentials have ``flat
directions,'' which is important for some questions in string
theory.  Also, in a forthcoming paper, we plan to extend our
methods to handle the class of manifolds with a Lie structure at
infinity \cite{ALN, MelroseScattering}, which is a class of
manifolds generalizing the class of manifolds with cylindrical
ends. Finally, our techniques and results may also be quite
relevant for problem arising in computational mathematics, more
precisely for obtaining fast algorithms for the finite-element
method (for solving the Dirichlet problem) on three dimensional
polyhedral domains, which is important in Science and Engineering.
See \cite{Arnold1, Arnold2} for an introduction to the finite
element method.

The reader is referred to \cite{ShubinBook, TaylorA}, or \cite{TaylorB} for
definitions and background material on pseudodifferential operators.
Note that in our paper we work exclusively with manifolds of bounded geometry.
The papers \cite{Shubin} and \cite{Strichartz} are a good introduction to some
basic results on the analysis on non-compact manifolds.  Throughout
the paper, a classical pseudodifferential operator $P$ will be called
{\em elliptic} if its principal symbols is invertible outside the zero
section.

Let us now briefly review the contents of each section (recall
that $M$ is a manifold with cylindrical ends). In Section
\ref{Sec.Operators} we introduce the algebra of operators
$\IPS{\infty}$ mentioned above and recall the classical
characterizations of Fredholm and compact operators in these
algebras. Section~\ref{Sec.spectral} deals with the same issues
for the algebra $\APS{\infty}$, which is a slight enlargement of
$\IPS{\infty}$, but has the advantage that it contains the
inverses of its elliptic, $L^2$-invertible elements.  We establish
several structure theorems for these algebras.  In
Section~\ref{Sec.potentials}, we introduce the double and single
layer potentials for manifolds with cylindrical ends and prove
that some of their basic properties continue to hold in this
setting.  In Section~\ref{Sec.parameter} we study boundary layer
potentials depending on a parameter on compact manifolds using a
method initially developed by G. Verchota in \cite{Ve}, and we
obtain estimates which are uniform in the parameter. These results
then allow us to establish the Fredholmness of the operators $S$
and $\pm\frac{1}{2}I+K$ discussed above. Finally, the last section
contains a proof of the Theorem~\ref{theorem.main}, which is a
statement about the well-posedness of the inhomogeneous Dirichlet
problem. This allows us to define and study the
Dirichlet-to-Neumann map in the same section.

\subsection{Acknowledgments} We are grateful to B. Ammann, T. Christiansen,
J. Gil, R. Lauter, G. Mendoza, and B. Monthubert for useful discussions.
We are also indebted to E. Schrohe who has sent us several of his
papers and answered some questions.

\section{Operators on manifolds with cylindrical ends\label{Sec.Operators}}

We begin by introducing the class of manifolds with cylindrical
ends (without boundary) and by reviewing some of the results on
the analysis on these manifolds that are needed in this paper.
Here we closely follow \cite{Guillemin, LauterMoroianu, Melrose46,
meaps, ScSc}. For simplicity, we shall usually drop the subscript
$M$ in the notation for the Laplacian $\Delta_{M}$ on $M$.

\subsection{Manifolds with cylindrical ends and the Laplace operator}
Let $M_1$ be a compact manifold with boundary $\pa M_1 \not =
\emptyset$. We assume that a metric $g$ is given on $M_1$ and that
$g_1$ is a product metric in a tubular neighborhood $V \cong \pa M_1
\times [0,1)$ of the boundary, namely
\begin{equation}
    g_1 = g_\pa + (dx)^2,
\end{equation}
where $x \in [0,1)$ is the second coordinate in $\pa M_1 \times [0,1)$
and $g_\pa$ is a metric on $\pa M_1$. Let
\begin{equation}
    M := M_1 \cup (\pa M_1 \times (-\infty, 0]), \quad \pa M_1
    \equiv \pa M_1 \times \{0\},
\end{equation}
be the union of $M_1$ and $\pa M_1 \times [0,\infty)$ along their
boundaries.  The above decomposition will be called a {\em
standard decomposition} of $M$.  The resulting manifold $M$ is called
a {\em manifold with cylindrical ends}.  Note that a manifold with
cylindrical ends is a complete, non-compact, Riemannian manifold
without boundary. Our class of manifolds with cylindrical ends is the
same as the one used in \cite{Shubincyl} in the framework of index theory.

Let $M = M_1 \cup (\pa M_1 \times (-\infty,0])$ be a manifold with
cylindrical ends. Let $g$ be the metric on $M$ and assume, as above,
that $g = g_\pa + (dx)^2$ on the cylindrical end $\pa M_1 \times
(-\infty,0]$, where $x \in (-\infty,0]$ and $g_\pa$ is a metric on the
boundary of $M_1$.  Let $d$ be the exterior derivative operator on $M$
so that $\Delta = \Delta_M = d^*d$ becomes the (scalar) Laplace
operator on $M$. Also, let $\Delta_{\pa M_1}$ be the Laplace operator
on $\pa M_1$, defined using the metric $g_\pa$. Then
\begin{equation}
    \Delta = \Delta_M = - \pa_x^2 + \Delta_{\pa M_1}
\end{equation}
on the cylindrical end $\pa M_1 \times (-\infty,0]$.

\subsection{Operators that are translation invariant in a neighborhood of
infinity} Let $M = M_1 \cup (\pa M_1 \times (-\infty,0])$ be a manifold
with cylindrical ends, as above, and let, for any $s \ge 0$,
\begin{equation}\label{eq.def.phis}
    \phi_s : \pa M_1 \times (-\infty, 0] \to \pa M_1 \times
    (-\infty, -s]
\end{equation}
be the isometry given by translation with $-s$ in the $x$-direction. If
$s<0$, then $\phi_s$ is defined as the inverse of $\phi_{-s}$. The special
form of the operator $\Delta$ obtained at the end of the previous
subsection suggests the following definition.

\begin{definition}\label{def.t.i}\
A continuous linear map $P : \CIc(M) \to \CI(M)$ will be called {\em
translation invariant in a neighborhood of infinity} if its Schwartz
kernel has support in
$$
    V_\epsilon := \{(x,y) \in M^2, \, \text{dist}(x,y) < \epsilon\},
$$
for some $\epsilon > 0$, and there exists $R > 0$ such that $P
\phi_s(f) = \phi_s P(f)$, for any $f \in \CIc(\pa M_1 \times (-\infty,
-R))$ and any $s > 0$.

We shall denote by $\IPS{m}$ the space of {\em order $m$, classical
pseudodifferential operators on $M$ that are translation invariant in a
neighborhood of infinity}.
\end{definition}

We have the following simple lemma.

\begin{lemma}\label{lemma.reg.b}\ Every $R \in \IPS{-n-1}$,
where $n$ is the dimension of $M$, induces a bounded operator on $L^2(M)$.
\end{lemma}

\begin{proof}\ The classical argument applies. Namely, $R$ is
defined by a continuous kernel $K$. The support condition on $K$
and the translation invariance at infinity then give
$$
    \int_M |K(x,y)|dx, \; \int_M |K(x,y)|dy \le C
$$
for some $C > 0$ that is independent of $x$ or $y$. This proves that
$R$ is bounded on $L^2(M)$, via Schur's lemma.
\end{proof}

We shall denote by $\mathcal D(T)$ the domain of a possibly unbounded
operator $T$.  Recall that an unbounded operator $T : \mathcal D(T)
\to X$ defined on a subset of a Banach space $Y$ and with values in
another Banach space $X$ is {\em Fredholm} if $T$ is {\em Fredholm} as
a bounded operator from its domain $\mathcal D(T)$ endowed with the
graph norm. Equivalently, $T$ is Fredholm if it is closed
and has finite dimensional kernel and cokernel. Also, $T$ is called
{\em invertible} if $T$ is invertible as an operator $\mathcal D(T)\to X$.
For all differential operators considered below, we shall consider the
minimal closed extension, that is, the closure of the operators with domain
compactly supported smooth functions.

For each nonnegative, even integer $m \in 2\NN$ we shall denote by $H^m(M)$
the domain of the operator $(I + \Delta)^{m/2}$ ($\Delta = \Delta_{M}$),
regarded as an unbounded operator on $L^2(M)$:

\begin{equation}\label{eq.def.H2m}
    H^{m}(M) := \mathcal D((I + \Delta)^{m/2}).
\end{equation}
We endow $H^m(M)$ with the norm
$$
    \|u \|_m = \| (I + \Delta)^{m/2} u \|_{L^2}
$$
(Below, we shall occasionally write $\|\;\cdot \;\|$ instead of
$\|\;\cdot \;\|_{L^2}$.) Note that $I + \Delta \ge I$, and hence
$$
    \|u \|_m \ge \|u\|.
$$
(Recall that $m \ge 0$.)

As usual, we shall denote by $H^{-m}(M)$ the dual of $H^m(M)$, via a
duality pairing that extends the pairing between functions and
distributions. We thus identify $H^{-m}(M)$ with a space of
distributions on $M$.

Let $\sigma_m(P) \in S^m(T^*M)/S^{m-1}(T^*M)$ be the principal symbol
of an operator $P \in \IPS{m}$. See \cite{ShubinBook, TaylorA}, or
\cite{TaylorB}.

\begin{lemma}\label{lemma.classical}\
Let $M$ be a manifold with cylindrical ends and $P \in \IPS{m}$ (so
$P$ is an order $m$ pseudodifferential operator that is translation
invariant in a neighborhood of infinity).
\begin{enumerate}[(i)]
\item\ For any $s, s'$, we have $\IPS{s}\IPS{s'} \subset \IPS{s+s'}$
and the principal symbol
$$
    \sigma_s : \IPS{s}/\IPS{s-1} \to S^s(T^*M)/S^{s-1}(T^*M)
$$
induces an isomorphism onto the subspace of symbols that are
translation invariant in a neighborhood of infinity.
\item\ Any $P \in \IPS{m}$ extends to a continuous map $P : H^{m'}(M)
\to H^{m'-m}(M)$, if $m, m' \in 2\ZZ$.
\end{enumerate}
\end{lemma}

\begin{proof}\ (i) follows from the analogous statement for
pseudodifferential operators on non-compact manifolds.

To prove (ii) when $m =0 $, we use the symbolic calculus,
Lemma~\ref{lemma.reg.b} and H\"ormander's trick. For $m'\ge m\ge 0$,
use the fact that $(I+\Delta)^k:L^2(M)\to H^{2k}(M)$ is an isometric
isomorphism and write
$$
    P = (I + \Delta)^i Q (I + \Delta)^j + R
$$
for suitable $Q,R \in \APS{0}$ and $i + j = m/2$. The other cases are
similar.
\end{proof}

Let us now recall a classical and well known construction.  Any
operator $P : \CIc(M) \to \CI(M)$ that is translation invariant in a
neighborhood of infinity will be properly supported (that is,
$P(\CIc(M)) \subset \CIc(M)$) and gives rise to a pseudodifferential
operator $\tilde{P}:\CIc(\pa M_1\times\RR)\to\CIc(\pa M_1\times\RR)$
by the formula
\begin{equation}\label{eq.def.tilde}
    \tilde{P} (f) = \phi_{-s} P \phi_s(f),
\end{equation}
where $\phi_s$ is the translation by $s$ on the cylinder $\pa M_1
\times \RR$ and $s$ is arbitrary, but large enough so that
$$
    \supp(P\phi_s(f))\, , \;\; \supp(\phi_s(f)) \subset \pa M_1
    \times (-\infty,0) \subset M.
$$
\vskip 0.08in

\begin{definition}\label{def.indicial}\
The operator $\tilde P$ will be called {\em the indicial operator
associated with $P$}. The resulting map
\begin{equation*} %\label{eq.def.Phi}
    \Phi: \IPS{\infty} \ni P \mapsto \tilde P \in \Psi^{\infty}(\pa
    M_1 \times \RR)
\end{equation*}
will be called the {\em indicial morphism}.
\end{definition}

Let us notice now that $\pa M_1 \times \RR$ is also a manifold with
cylindrical ends. The partially defined action of $\RR$ on the ends of $M$
extends to a global action of $\RR$ on $\pa M_1 \times \RR$.
Let $T\in\Psi^{-\infty}_{\rm inv}(M_1 \times \RR)^{\RR}$ and $\eta$ be a
smooth function on $\RR \times \pa M_1$ with support in $(-\infty,-1)\times
\pa M_1$, equal to $1$ in a neighborhood of infinity. Then
\begin{equation}\label{eq.def.sec}
    s_0(T) : = \eta T \eta
\end{equation}
defines an operator in $\IPS{-\infty}$.

Let us denote by $\Psi^{\infty}_{\rm inv}(\pa M_1 \times \RR)^{\RR}$
the operators in $\Psi^{\infty}_{\rm inv}(\pa M_1 \times \RR)$ that
are translation invariant with respect to the natural action of $\RR$
on $\pa M_1 \times \RR$.

\begin{lemma}\label{prop.ran.Phi}\ Let $s_0$ be as in Equation
\eqref{eq.def.sec}.  Then $\Phi(s_0(T)) = T$ for all $T \in
\Psi^{\infty}_{\rm inv}(\pa M_1 \times \RR)^{\RR}$.  In particular,
the range of the indicial morphism $\Phi$ of Definition
\ref{def.indicial} is $\Psi^{\infty}_{\rm inv}(\pa M_1 \times
\RR)^{\RR}$.
\end{lemma}

\begin{proof}\ This is a direct consequence of the definition.
\end{proof}

In order to deal with operators acting on weighted Sobolev spaces,
we shall need the following lemma. (See also \cite{meaps}.)

\begin{lemma}\label{lemma.simple}\ Let $P,P_1 \in \IPS{\infty}$ be
arbitrary and $\rho : M \to [1,\infty)$ be a smooth function such that
$\rho(y,x)=x$ on a neighborhood of infinity in $\pa M_1 \times (-\infty,0]$.
Then
\begin{enumerate}[(i)]
    \item\  $\tilde Q = \tilde P \tilde {P_1}$, if $Q = PP_1$.
    \item\  $\adj_\rho(P) := [\rho, P ] \in \IPS{\infty}$.
\end{enumerate}
\end{lemma}

\begin{proof}\ The relation (i) follows by chasing definitions. To prove
(ii), we can assume that $M = X \times \RR$. Let $\phi_s$, $s\in\RR$,
be translation by $s$ along $\RR$. We can assume that $P$ is translation
invariant, in the sense that $\phi_s^*(P) = P$, for any $s>0$. Then
\begin{equation}
    \phi_s^*([x,P]) = [\phi_s^*(x), \phi_s^*(P)] = [ x + s, P ] =
    [ x,P ].
\end{equation}
Thus $[x,P]$ is also $\RR$-invariant.

In general, $\rho = x$ in a neighborhood of the infinity, so the
result follows.
\end{proof}

The properties of the indicial operators $\tilde P$ are conveniently
studied in terms of {\em indicial families}. Indeed, by considering
the Fourier transform in the $\RR$ variable, we obtain by Plancherel's
theorem an isometric bijection (that is, a unitary operator)
defined, using local coordinates $y$ on $\pa M_1$, by

\begin{equation}\label{fourier}
    \mathcal{F}:L^2(\pa M_1\times\RR)\to L^2(\pa M_1\times\RR),\quad
    \mathcal{F}(f)(y,\tau):=
    \frac{1}{\sqrt{2\pi}}\int_\RR e^{-\imath\tau x}f(y,x)\,dx.
\end{equation}

\noindent Hereafter, $\imath:=\sqrt{-1}$.

Because $\tilde P$ is translation invariant with respect to
the action of $\RR$, the resulting operator $P_1:=\mathcal F\tilde{P}
\mathcal{F}^{-1}$ will commute with the multiplication operators in $\tau$,
and hence it is a decomposable operator, in the sense that there exist
(possibly unbounded) operators $\hat P(\tau)$ acting on $\CI(\pa M_1)
\subset L^2(\pa M_1)$ such that
$$
    (P_1 f) (\tau) = \hat P(\tau) f(\tau), \quad f(\tau) =
    f(\cdot,\tau) \in \CI(\pa M_1).
$$
In other words,

\begin{equation}\label{multiplier}
    \bigl[(\mathcal{F}\tilde{P}\mathcal{F}^{-1}f)\bigr](\tau)
    = \hat{P}(\tau)f(\tau).
\end{equation}

\noindent Using local coordinates, it is not hard to see that the operators
$\hat P(\tau)$ are classical pseudodifferential operators and that the map
$\tau\mapsto\hat P(\tau)f$ is $\CI$ for any $f \in \CI(\pa M_1)$.

One also has $\tilde P(e^{\imath\tau x} g)= e^{\imath \tau x}\hat{P}(\tau)g$,
for any $g \in L^2(\pa M_1)$. Let $K_{\tilde P}$ be the distribution kernel
of $\tilde P$. Then
\begin{equation}
    K_{\tilde P} (x_1, x_2, y_1, y_2) = k_{\tilde P}(x_1 - x_2, y_1, y_2),
\end{equation}
for some distribution $k_{\tilde P}$ on $\RR \times (\pa M_1)^2$.
This allows us to write the distribution kernel of $\hat P(\tau)$ as
\begin{equation}
    K_{\hat P(\tau)}( y_1, y_2) =
    \int_\RR k( x, y_1, y_2)\,e^{ - \imath t\tau}\,dx.
\end{equation}
Let $Q = [\rho, P]$. Then
\begin{equation}\label{eq.ind.com}
    k_{\tilde Q} = \imath \frac{\pa}{\pa \tau} \hat P(\tau).
\end{equation}
See \cite{LauterMoroianu, Melrose46, meaps} and the references
therein.

\section{A spectrally invariant algebra\label{Sec.spectral}}

A serious drawback of the algebra $\IPS{\infty}$ is that it is not
``spectrally invariant,'' in the sense that the inverse of an elliptic operator
$P\in\IPS{\infty}$ that is invertible on $L^2$ is not necessarily in
this algebra (Definition \ref{def.spec.inv} below). In this section we
slightly enlarge the algebra $\IPS{\infty}$ so that it becomes
spectrally invariant. This will lead us to an algebra of operators that
are ``almost translation invariant in a neighborhood of infinity.''

\subsection{Operators that are almost translation invariant in a
neighborhood of infinity} We begin by introducing another algebra of
pseudodifferential operators that will be indispensable also later
on. Let $\rho$ be a smooth function as in Lemma \ref{lemma.simple}.
Recall that $\adj_\rho(T):=[\rho,T]$. Assume $T:\CIc(M) \to \CI(M)$
to be a linear map with the property that
$$
    \adj^k_\rho(T) := [\rho, [\rho, \ldots, [\rho, T] \ldots ]]
$$
extends to a continuous map $\adj^k_\rho(T): H^{-m}(M) \to H^{m}(M)$,
for any $m \in 2\ZZ_+$. Let $\|T\|_{k,m}$ denote the norm of the
resulting operator $\adj^k_\rho(T)$. Recall the section $s_0$ defined
in Equation \eqref{eq.def.sec}.

We define $\APS{-\infty}$ to be the closure of $\IPS{-\infty}$ with
respect to the countable family of semi-norms
\begin{equation}\label{eq.def.fam.s}
    T \to \|T \|_{k,m}\,, \quad \text{and} \quad T \to
    \|\rho^l(T - s_0(\Phi(T)))\rho^{l} \|_{0,m}.
\end{equation}
where $k, m/2,l \in \ZZ_+$. Then $\APS{-\infty}$ is a
Fr\'echet algebra (that is, a Fr\'echet space endowed with an algebra
structure such that the multiplication is continuous).

Finally, we define
\begin{equation}
    \APS{m} := \IPS{m} + \APS{-\infty}.
\end{equation}
An element $P \in \APS{m}$ will be called {\em almost translation
invariant in a neighborhood of infinity}.

It is interesting to observe now that we can introduce dependence
on $\rho$ at infinity (thus obtaining variants of Melrose's
b-calculus, see \cite{meaps} and \cite{LMN}). This is done by
noticing that for any $P \in \IPS{m}$ and any $N \in \NN$ there
exists a bounded operator $R_N : H^{-k}(M) \to H^{k}(M)$, where
$2k \le m - N$, such that
\begin{equation}
    \rho^{-a} P \rho^a - \sum_{j = 0}^{N - 1} (-1)^j \rho^{-j}
    \binom{a}{j}\adj_\rho^j(P) = \rho^{N/2}R_N \rho^{N/2}.
\end{equation}
(Above, $\binom{a}{j}=a(a-1) \ldots (a-j+1)/j!$ stand for the usual
``binomial'' coefficients.)

We now define the fractionary Sobolev spaces. Let $s \ge 0$ and choose
$P_s \in \APS{s}$ to be elliptic and to satisfy $P_s \ge 1$. We shall
denote by $H^s(M)$ the domain of (the closure of) $P_s$, regarded as
an unbounded operator on $L^2(M)$:
\begin{equation}\label{eq.def.Hs}
    H^s(M) := \mathcal D(\overline{P}_s).
\end{equation}
This definition is independent of our particular choice of $P_s$
because, if $P_s'$ is another such selection, we can choose $Q \in
\APS{0}$ and $R \in \APS{-\infty}$ such that
\begin{equation}
    P_s' = Q P_s + R.
\end{equation}
Thus, if $\xi \in \mathcal D(P_s)$, then there exists a sequence
$\xi_n \in \CIc(M)$, $\xi_n \to \xi$ in $L^2(M)$, such that $P_s\xi_n$
converges in $L^2(M)$.  But then $P_s'(\xi_n) = Q(P_s\xi_n) + R\xi_n$
also converges, because $Q$ and $R$ are continuous. See also
\cite{LMN}.

We endow $H^s(M)$ with the norm $\| f \|_s := \|P_s
f\|_{L^2(M)}$. (Using a quantization map from symbols to
pseudodifferential operators, we can assume that $\|f\|_s$ depends
analytically on $s$.)  For $s < 0$, $H^s(M)$ is the dual of
$H^{-s}(M)$, regarded as a space of distributions on $M$.  The
subspace $\CIc(M) \subset H^s(M)$ is dense. See \cite{ALNV} for more
results on Sobolev spaces on manifolds with a Lie structure at
infinity, a class of manifolds that includes the class of manifolds
with cylindrical ends. For example, $H^s(M)$ can be identified with
the domain of $(I+\Delta)^{s/2}$.

We shall also consider {\em weighted Sobolev spaces} as follows. Let
$\rho : M \to [1, \infty)$, $\epsilon > 0$, be a smooth functions such
that $\rho(y,x) = x$, for $(y,x) \in \pa M_1 \times (-\infty,-R]$ with
$R $ large enough, as before. Then we shall denote by $\rho^a H^s(M)$ the space
of distributions of the form $\rho^a u$, with $u \in H^s(M)$. We endow
$\rho^a H^s(M)$ with the norm
$$
    \|f\|_{s,a} := \|\rho^{-a}f\|_s.
$$

We have then the following classical results about almost translation
invariant pseudodifferential operators on the manifold with
cylindrical ends $M$ \cite{MeMe}.  (See \cite{Kondratiev, LMN,
LockhartOwen, MeN, Schrohe, ScSc, Schulze}.)  These results generalize
the corresponding even more classical results on pseudodifferential
operators on compact manifolds.

\begin{theorem}\label{theorem.classical2}\
Let $M$ be a manifold with cylindrical ends and $P \in \APS{m}$
(so $P$ is an order $m$ pseudodifferential operator that is almost
translation invariant in a neighborhood of infinity). Also, let $\rho>0$,
$\rho(y,x) = x$ on a neighborhood of infinity in $\pa M_1 \times
(-\infty, 0]$. Let $s,a \in \RR$ be arbitrary, but fixed. Then:
\begin{enumerate}[(i)]
\item\  $P$ extends to a continuous operator $P : \rho^a H^s(M) \to
\rho^{a} H^{s-m}(M)$.
\item\  $P : \rho^a H^s(M) \to \rho^{a'} H^{s-m'}(M)$ is compact for any
$a' < a$ and $m' > m$.
\item\  $P : \rho^a H^s(M) \to \rho^a H^{s-m}(M)$ is compact
$\Leftrightarrow \; \sigma_m(P) = 0$ and $\tilde P = 0$.
\item\  $P : \rho^a H^s(M) \to \rho^a H^{s-m}(M)$ is Fredholm
$\Leftrightarrow \; \sigma_m(P)$ is invertible and the operator
$\tilde P: H^s(\pa M_1 \times\RR) \to H^{s-m}(\pa M_1 \times \RR)$
is an isomorphism.
\end{enumerate}
\end{theorem}

\begin{proof}\ This theorem follows for example from the results in
\cite{meaps}, or the older preprint \cite{MeMe}.
\end{proof}

A far reaching program for generalizing the above result to other
classes of non-compact manifolds is contained in Melrose's ``small
red book'' \cite{MelroseScattering}. See also \cite{Nistor-prden}.
Also, see \cite{SchSe} for an extension of the above results to
$L^p$--spaces, and \cite{CoSchSe} for some applications to
non-linear evolution equations.

As a consequence, we obtain the following result.

\begin{corollary}\label{cor.selfadj}\ Let $P \in \APS{m}$, $m > 0$, be
elliptic. If $P$ is symmetric on $\CIc(M)$, then it is essentially
self-adjoint.
\end{corollary}

\begin{proof}\
We replace $P$ by its closure first.
We want to prove that $P \pm \imath I$ is invertible.
Denote the inner product on $L^2(M)$ by
$\<\,\cdot\,,\,\cdot\,\>$. Then $\<(P \pm \imath I)\xi, \xi \> =
\|P\xi\|^2 + \|\xi\|^2$, for any $\xi$ in the domain of $P$, and hence
$P \pm \imath I$ is injective and has closed range.

Let us prove that the range of $P \pm \imath I$ is dense. We deal only
with $P + \imath I$, because the other case is completely
similar. Assume the range of $P \pm \imath I$ is not dense,
then there exists $\eta \in L^2(M)$ such that
$$
    \<(P + \imath I)\xi,\eta\> = 0
$$
for all $\xi \in \CIc(M)$. Then $(P-\imath I)\eta = 0$ in the sense of
distributions. Select $Q\in\APS{-m}$ such that $Q(P-\imath I)=I-R$, where
$R \in \APS{-\infty}$. Then $\eta = R\eta$. Choose
$\eta_n \in \CIc(M)$, $\eta_n \to \eta$ in $L^2(M)$. By the definition
of $\APS{-\infty}$, we can find operators $R_n \in \IPS{-\infty}$ such that
$$
    \|R-R_n\|_{0,m'} := \|(I+\Delta)^{m'/2} (R-R_n)
    (I+\Delta)^{m'/2}\|\to 0,
$$
for $m' \ge m$.  Then $\xi_n := R_n\eta_n \to \eta$, as well, and
$\xi_n \in \CIc(M)$. Moreover,
$$
    (P - \imath I)\xi_n = (P - \imath I) R_n \eta_n \to (P -
    \imath I)R \eta,
$$
because the operators $(P - \imath I) R_n$ are bounded and converge in
norm to $(P - \imath I) R \in \APS{-\infty}$. This proves that $\eta$
is in the domain of the closure of $P$, which is a contradiction,
since we have already seen that $P - \imath I$ is injective.
\end{proof}

We now investigate the structure of the ideals of the algebras
$\IPS{-\infty}$ and, most important, $\APS{-\infty}$.

\begin{lemma}\label{lemma.quot}\ The range of the map
\begin{equation}\label{eq.def.Phi}
    \Phi : \APS{-\infty} \ni P \mapsto \tilde P \in
    \Psi^{-\infty}(\pa M_1 \times \RR)
\end{equation}
identifies with $\mathcal S(\RR \times (\pa M_1)^2)$, via the map
$\chi$ that sends the kernel $K(t_1,t_2,y_1,y_2) \in \CI(\RR^2 \times
(\pa M_1)^2)$ of $\tilde P$ to the function $k(t,y_1,y_2) = K(t, 0,
y_1, y_2) \in \mathcal S(\RR \times (\pa M_1)^2)$.

In particular, $\Phi(\APS{-\infty}) = \Psi^{-\infty}_{\rm ai}(\pa M_1
\times \RR)^{\RR}$.
\end{lemma}

\begin{proof}\ The indicial map
$$
    \Phi : \IPS{-\infty} \to \Psi^{-\infty}_{\rm inv}(\pa M_1 \times \RR)
$$
of Definition \ref{def.indicial} is by definition continuous. It is
also surjective by Proposition \ref{prop.ran.Phi}. It has a canonical
continuous section $s_0$, which associates to $T \in \Psi_{\rm
inv}^{\infty}(M_1 \times \RR)^\RR$ the operator $s_0(T) := \eta T
\eta$, where $\eta$ is a smooth function on $\RR \times \pa M_1$ and
with support in $(-\infty, -1) \times \pa M_1$, and equal to $1$ in a
neighborhood of infinity (cf. Equation \eqref{eq.def.sec}).

Moreover, $s_0$ sends properly supported operators to $\IPS{-\infty}$.
This shows that
$$
    \IPS{-\infty}\cong \ker(\Phi) \oplus s_0(\Psi_{\rm inv}^{\infty}
    (M_1 \times \RR))^\RR,
$$
as Fr\'echet spaces. We also see that the quotient seminorms defined
by the seminorms of Equation \eqref{eq.def.fam.s} on the range of
$\Phi$ are the same as the seminorms defining the topology on
$\mathcal S(\RR \times (\pa M_1)^2)$.  Since $\mathcal S(\RR \times
(\pa M_1)^2)$ is the closure of $\chi(\Psi^{-\infty}_{\rm inv}(M_1
\times \RR))$, the result follows.
\end{proof}

The same proof as above also gives the following result.

\begin{corollary}\label{cor.desc}\
Let $Y$ be a compact, smooth manifold without boundary. Then the
algebra $\Psi^{-\infty}_{ai}(Y \times \RR)^{\RR}$ is the space of
operators $T$ on $L^2(Y \times \RR)$ such that $(I + \Delta_{Y \times
\RR})^m \adj_\rho^k(T) (I + \Delta_{Y \times \RR})^m$ is bounded and
$\RR$--invariant for any $m, k \in \ZZ_+$.  The resulting family of
seminorms is the family of seminorms of Equation \eqref{eq.def.fam.s}
defining the topology on $\Psi^{-\infty}_{ai}(Y \times \RR)^{\RR}$.
\end{corollary}

Let $\mfk I$ be the kernel of the map $\Phi : \APS{-\infty} \to
\Psi^{-\infty}(\pa M_1 \times \RR)$ of Definition
\eqref{eq.def.Phi}. We also have the following description of $\mfk I$
that is similar in spirit to Corollary \ref{cor.desc}.

\begin{lemma}\label{lemma.ideal}\
The ideal $\mfk I$ is the space of all operators $T$ on $L^2(M)$ such
that $(I + \Delta)^m \rho^l T \rho^l (I + \Delta)^m$ is bounded for
any $m, l \in \ZZ_+$. The resulting family of seminorms is the family
of seminorms of Equation \eqref{eq.def.fam.s} defining the topology on
$\mfk I$.
\end{lemma}

\begin{proof}\ It is clear from the definition that
$$
    T \to \|(I + \Delta)^m \rho^l T \rho^l (I + \Delta)^m\|
$$
is one of the seminorms of Equation \eqref{eq.def.fam.s}, namely
$\|\;\cdot \;\|_{l,m}$.

Conversely, let $T$ be an operator on $L^2(M)$ such that for each
$m,l\in\ZZ_+$ the operator $(I + \Delta)^m \rho^l T \rho^l (I +
\Delta)^m$ is bounded.  The family of seminorms
$T\to\|\rho^l(I+\Delta)^m T(I+\Delta)^m\rho^l\|$ is equivalent to the
family $\|\;\cdot\;\|_{l,m}$. We shall use this family instead.

The Schwartz kernel of $T$ is $K_T(x,y) = \< T \delta_y, \delta_x\>$
and it satisfies
\begin{equation}\label{eq.eq.l}
    \rho^l(x) \rho^l(y) |K_T(x,y)| \le C^2\|\rho^l(I + \Delta)^m T
    (I + \Delta)^m\rho^l\|
\end{equation}
where $C \ge \|\delta_x\|_{-m}$, uniformly in $x \in M$, for some $m >
n/2$. (We have used here the Sobolev embedding theorem for manifolds
with cylindrical ends \cite{ALNV}.)

We shall prove now that $T$ is in the closure of $\ker \Phi \subset
\IPS{-\infty}$ (recall that $\Phi(T) = \tilde T$ is the indicial map).
Let $\alpha_n = 1 - \phi_{n}(\eta) \in \CIc(M)$, where $\phi_n$ is
translation by $-n$ on the cylindrical end, and $\eta \in \CI(M)$ is
equal to $1$ in a neighborhood of infinity and is supported on $\pa
M_1 \times (-\infty, 0]$, if
$$
    M = M_1 \cup \pa M_1 \times (-\infty, 0]
$$
is a standard decomposition of $M$.

We have that $T_n := \alpha_n T \alpha_n$ has the compactly supported
Schwartz kernel
$$
    K_{T_n}(x,y) = \alpha_n(x) K_T(x,y) \alpha_n(y).
$$
Taking $l > 1$ in the Equation \eqref{eq.eq.l}, we see using Shur's lemma
(as in the proof of Lemma \ref{lemma.reg.b}) that $\|T_n - T\| \to 0$
(the norm here is that of bounded operators on $L^2(M)$). The proof
that $\|T_n - T\|_{l,m} \to 0$ for $l > 0$ or $m > 0$ is completely
similar.
\end{proof}

Let $M = M_1 \cup (\pa M_1 \times (-\infty, 0])$ be a standard
decomposition of $M$. Consider a diffeomorphism $\psi$ from $M$ to the
interior of $M_1$ that coincides with $(y,t) \mapsto (y,-t^{-1})$ in a
neighborhood of infinity.

\begin{corollary}\label{cor.ideal}\ The diffeomorphism $\psi$ above
identifies $\mfk I$ with $\CI_0(M_1^2)$, that is, the space of smooth
functions on $M_1^2$ that vanish to infinite order at the boundary
$\pa (M_1^2) = (\pa M_1 \times M_1) \cup (M_1 \times \pa M_1)$.
\end{corollary}

\begin{proof}\
This follows right away from the proof of Lemma \ref{lemma.ideal}.
\end{proof}

To formulate the following results, it is convenient to use the
following classical concept (see \cite{Schrohe2}, for example).

\begin{definition}\label{def.spec.inv}\
Let $\mathcal A$ be an algebra of bounded operators on some Hilbert
space $\mathcal H$. We say that $\mathcal A$ is {\em spectrally
invariant} if, and only if, $(I + T)^{-1} \in I +\mathcal A$, for any
$T \in \mathcal A$ such that $I + T$ is invertible as an operator on
$\mathcal H$.
\end{definition}

\begin{lemma}\label{lemma.spec.inv}\
The algebras $\mfk I$ and $\Psi_{\rm ai}^{-\infty} (\pa M_1 \times
\RR)^{\RR}$ are spectrally invariant.
\end{lemma}

\begin{proof}\ Both are well known results (see \cite{LMN} or
\cite{ScSc} and the references therein). An easy proof is obtained
using Lemma \ref{lemma.ideal} or, respectively, Corollary \ref{cor.desc}.
\end{proof}

The property of being spectrally invariant is preserved under
extensions of algebras (see \cite{LMN}). Using this twice, we obtain
the following result.

\begin{corollary}\label{cor.sp.i}\ The algebras $\APS{-\infty}$ and
$\APS{0}$ are spectrally invariant.
\end{corollary}

A proof of this corollary is also contained in the following theorem,
which is the main result of this section. It states
that $\APS{\infty}$ is, in a certain sense, also spectrally invariant,
its proof does not rely on the above corollary.

\begin{theorem}\label{theorem.sp.i}\ Let $T \in \APS{m}$,
$m \ge 0$, be such that $T$ is invertible as a (possibly unbounded)
operator on $L^2(M)$. If $m > 0$, we assume also that $T$ is
elliptic. Then $T^{-1} \in \APS{-m}$.
\end{theorem}

\begin{proof}\ Note that for $m=0$, it is a consequence of the
invertibility that $T$ must again be elliptic, as in the case
$m >0$.

Let $Q_1$ be a parametrix of $T$, namely, $Q_1 \in \APS{-m}$ and
$$
    Q_1P - I\, , \, \; PQ_1 - I \in \APS{-\infty}.
$$
Let $\xi$ be a distribution such that $\xi\, , \,\, T\xi \in L^2(M)$. Then
$$
    \xi = Q_1 (P\xi) - (Q_1P - I) \xi \in H^m(M).
$$
This shows that the maximal domain of $T$ is $H^m(M)$. Since $T$ is invertible,
the graph topology on the domain of $T$ coincides with the topology
of $H^m(M)$. It follows then that $T : H^m(M) \to L^2(M)$ is Fredholm
(in fact, even invertible) and hence $\hat T(\tau)$ is invertible
for any $\tau \in \RR$.

Let $R_1 \in \APS{-\infty}$ be such that $R_1$ and $\hat R_1(\tau)$ are
injective, for any $\tau \in \RR$. (This is possible because $L^2(M)$
has a countable orthonormal basis.)  Then $Q_2 := Q_1^*Q_1 + R_1^*R_1$ is
a parametrix of $TT^*$ such that $Q_2 : L^2(M) \to H^{2m}(M)$ is an
isomorphism. Let $R_2 := TT^*Q_2 - I \in \APS{-\infty}$. By construction,
$I + R_2$ is invertible on $L^2(M)$ and $I + \hat R_2(\tau)$ are invertible
on $L^2(\pa M_1)$ for any $\tau \in \RR$.

By Lemma \ref{lemma.spec.inv} applied to $ I + \hat R_2(\tau)$ and the
algebra $\Psi_{\rm ai}^{-\infty}(\pa M_1 \times \RR)^{\RR}$, we can find
$R_3 \in \APS{-\infty}$ such that $(I + R_2)(I + R_3) - I \in \mfk I$.
(We can take $R_3 = s_0[(I + \tilde R_2)^{-1} - I]$. We
can also assume that $I + R_3$ is injective, by replacing $I + R_3$ with
$$
    (I + R_2)^*[(I + R_3)^*(I + R_3) + R_4^*R_4],
$$
where $R_4 \in \mfk I$ is injective.

We are now ready to complete our proof. The operator
$R_5:=TT^*Q_2(I+R_3)-I\in \mfk I$ is such that $I + R_5$ is injective.
It follows that $I + R_5 =TT^*Q_2(I + R_3)$ is Fredholm of index zero
and, hence, invertible on $L^2(M)$. Using again Lemma \ref{lemma.spec.inv},
we obtain that there exists $R_6 \in \mfk I$ such that
$(I + R_5)(I +R_6) = I$. Thus,
$$
    TT^* Q_2(I + R_3)(I + R_6) = I.
$$
This means that $P : = T^*Q_2(I + R_3)(I + R_6)$ is a right inverse to $T$.
We can prove in exactly the same way that $T$ has a left inverse in
$\APS{-m}$ and, hence, that it is invertible in $\APS{\infty}$.
\end{proof}

The above theorem applied to $T = I + \Delta$ gives the following result.

\begin{corollary}\label{cor.Delta-1}\
Let $M$ be a manifold with cylindrical ends and $\Delta = \Delta_M$ be the
Laplace operator on $M$. Then $(I + \Delta)^{-1} \in \APS{-2}$.
\end{corollary}

\subsection{Perturbation by potentials}
We shall need also a further extension of the above corollary.  To
state it, recall that an operator $L$, mapping $L^2_{\rm loc}$ into
distributions, is said to have the {\it unique continuation property}
if
$$
    Lu=0\quad\&\,\,\mbox{ $u$ vanishes in an open set}
    \Longrightarrow\,u=0\mbox{ on }M.
$$

\begin{proposition}\label{prop.L+V}\ Let $L \in \APS{m}$ be nonnegative
(that is, $L \ge 0$) and satisfy the unique continuation
property. Also, let $V \in \CI(M) \cap \APS{0}$ (that is, $V$ is
translation invariant in a neighborhood of infinity), $V \ge 0$, such
that $V$ is strictly positive on some open subset of $M$. Then, if $L
+ V : H^m(M) \to L^2(M)$ is Fredholm, it is also invertible.
\end{proposition}

\begin{proof}\
The assumptions $L \ge 0$ and $V \ge 0$ imply $L + V \ge 0$, as
well. Assume by contradiction that $L + V: H^m(M) \to L^2(M)$ is
Fredholm but not invertible. Then $L + V$ is Fredholm as an unbounded
operator on $L^2(M)$ and is not invertible. This shows that $0$ must be
an eigenvalue of $L + V$.

Let $u \not = 0,$ $u \in L^2(M)$ be an associated eigenvector:
$$
    (L + V)u = 0.
$$
Then
$$
    \<Lu,u\> + \<Vu,u\> = 0,
$$
where $\<\cdot,\cdot\>$ is the inner product on $L^2(M)$.

Since $\<Lu,u\> \ge 0$ and $\<Vu,u\> \ge 0$, we must have both $\|L^{1/2}u
\|^2 = \<Lu,u\> = 0$ and $\|V^{1/2}u\|^2 = \<Vu,u\> = 0$. Thus $Lu = 0$
and $Vu = 0$. The second relation gives that $u$ vanishes on some open
subset of $M$. Since $L$ has the unique continuation property, $u$
must vanish identically. This contradicts the original assumptions and
the proof is now complete.
\end{proof}

\vskip 0.08in
\noindent{\bf Example.} If $T\in \APS{k}$ has the unique continuation
property then $L:=T^*T$ satisfies the hypotheses of the above
proposition (with $m=2k$). In particular, this is the case for $\Delta
= \Delta_M=d^*d$, since the kernel of $d=d_M$ consists of only locally
constant functions.
\vskip 0.08in

The following theorem is crucial for our approach to extending the
method of layer potentials to manifolds with cylindrical ends.

\begin{theorem}\label{theorem.Delta+V}\
Let $M$ be a manifold with cylindrical ends and $V \ge 0$ be a smooth
function on $M$ that is translation invariant in a neighborhood of
infinity and does not vanish at infinity. Denote by $\Delta= \Delta_M$ the
Laplace operator on $M$.  Then $\Delta + V$ is invertible as an
unbounded operator on $L^2(M)$ and $(\Delta + V)^{-1} \in \APS{-2}$.
\end{theorem}

\begin{proof}\ For starters, $\Delta$ is non-negative ($\Delta \ge 0$) and
has the unique continuation property (cf. the previous example).
Since the potential $V$ is non-negative, as well as strictly positive on
some non-empty open set, our result will follow from
Proposition~\ref{prop.L+V} as soon as we show
that $\Delta + V : H^2(M) \to L^2(M)$ is Fredholm.

Since $\Delta$ is elliptic, $P := \Delta + V: H^2(M) \to L^2(M)$
will be Fredholm if, and only if, $\tilde P$ is invertible. In turn,
to show that $\tilde P$ is invertible it suffices to prove the norm
of the inverse of $\hat P(\tau) : H^2(\pa M_1) \to L^2(\pa M_1)$ is
bounded uniformly in $\tau\in \RR$.

More specifically, let $V_\infty \in\CI(\pa M_1)$ be the limit at
infinity of the function $V$. (This limit exists because we assumed
$V$ to be translation invariant in a neighborhood of infinity.)
Denote $\Delta = \Delta_{\pa M_1}$, to simplify notation in what follows.
By definition, we have
$$
    \hat P(\tau) = \Delta + \tau^2 + V_\infty.
$$
Since $V_\infty + \tau^2 \ge 0$ and does not vanish identically
for any $\tau \in \RR$, by assumption, we obtain as in \cite{MTJFA}
that $\hat P(\tau)$ is indeed invertible for any $\tau \in \RR$.
(One can also justify this using the methods used to prove
Proposition~\ref{prop.L+V}.)

Let $\mathcal{L}(X,Y)$ denote the normed space of all linear
bounded operators between two Banach spaces $X$, $Y$.

The invertibility of $\Delta + V_\infty$ implies that
$\Delta_{\pa M_1} + V_\infty \ge cI$, for some $c > 0$. The functional
calculus gives that $(\Delta + \tau^2 + V_\infty)^2 \ge c^2I$ and that
\begin{equation*}
    (\Delta + \tau^2 + V_\infty)^2 \ge (\Delta + V_\infty)^2 =
    \frac{1}{2} \Delta^2 + \frac{1}{2} (\Delta + V_\infty)^2 - V_\infty^2
    \ge \frac{1}{2} \Delta^2 - \|V_\infty\|_\infty^2\,.
\end{equation*}
Consequently,
\begin{equation}\label{eq.ineq.Delta}
    (\Delta + \tau^2 + V_\infty)^2 \ge \frac{\epsilon}{2}
    (\Delta^2-2\|V_\infty\|_\infty)
    + (1-\epsilon)c^2 \ge 2C^2(\Delta^2 + 1) \ge C^2(\Delta + 1)^2,
\end{equation}
if $\epsilon>0$ and $C>0$ are small enough. In particular, we obtain from
Equation~\eqref{eq.ineq.Delta} that
$$
    \|(\Delta + \tau^2 + V_\infty)(\Delta+1)^{-1}\| \ge C,
$$
and, ultimately,
\begin{multline*}
    \| (\Delta_{\pa M_1} + \tau^2 + V_\infty)^{-1}
    \|_{\mathcal{L}(L^2(\pa M_1),H^2(\pa M_1))} \\ = \|
    (\Delta_{\pa M_1} + 1) (\Delta_{\pa M_1} + \tau^2 +
    V_\infty)^{-1} \|_{\mathcal{L}(L^2(\pa M_1),L^2(\pa M_1))} \le
    C^{-1},
\end{multline*}
for any $\tau \in \RR$. This completes the proof of our theorem.
\end{proof}

Let us mention that in the proof of the above theorem we
used an ad-hoc argument to prove a result that holds in much greater
generality.  Namely, assume that $P$ is elliptic of order $m$. Then
there exists $R >0$ such that $\hat P(\tau)$ is invertible as a map
$H^m(\pa M_1) \to L^2(\pa M_1)$, for any $|\tau| > R$. Moreover, $\hat
P(\tau)^{-1}$ depends continuously on $\tau$ on its domain of definition.
In particular, if $P$ is elliptic of order $m > 0$ and $\hat P(\tau)$
is invertible  for any $\tau$, then $\| \hat P (\tau)^{-1} \|$ is
uniformly bounded as a map $L^2(\pa M_1) \to H^{m}(\pa M_1)$.
See \cite{Shubin}, especially Theorem 9.2, for details.

\subsection{Products} We shall need also the following product
decomposition result for the ideal of regularizing, almost invariant
pseudodifferential operators.

First, let us observe that if $M$ is a manifold with cylindrical ends
and $X$ is a smooth, compact, Riemannian manifold without boundary,
then $M \times X$ is also a manifold with cylindrical ends.

For any Fr\'echet algebra $A$, we shall denote by $\CI(X^2,A)$ the
space of smooth functions on $X \times X$ and values in $A$, with the
induced topology and the product:
\begin{equation}
    f \star g (x,x'') = \int_X f(x,x') g(x',x'')dx' \,,
\end{equation}
the integration being with respect to the volume element obtained from
the Riemannian metric on $X$. For example, $\Psi^{-\infty}(X)\cong
\CI(X^2,\CC)$.

\begin{theorem}\label{theorem.prod}\
Let $M$ be a manifold with cylindrical ends and $X$ be a smooth,
compact, Riemannian manifold without boundary. Then $\APSX{-\infty}$
is isomorphic to $\CI(X^2,\APS{-\infty})$.
\end{theorem}

\begin{proof}\ Let us denote by $\mathcal S(\RR, V)$ the Schwartz space
of rapidly decreasing smooth functions on $\RR$ with values in a
Fr\'echet space $V$. Also, let $\CI_0$ denote the space of smooth functions
on a manifold with boundary that vanish to infinite order at the boundary,
as in the statement of Corollary \ref{cor.ideal}.

The statement of the theorem follows from Lemma \ref{lemma.quot}, Corollary
\ref{cor.ideal}, and the relations
\begin{equation}
\begin{gathered}
    \mathcal S(\RR, \CI( (\pa M_1 \times X)^2 )) \simeq \CI( X^2,
    \mathcal S(\RR, \CI( (\pa M_1)^2 )))\,, \quad \text{and} \\
    \CI_0( (\pa M_1 \times X)^2 ) \simeq \CI( X^2, \CI_0( (\pa M_1)^2))\,.
\end{gathered}
\end{equation}
\end{proof}

\section{Boundary layer potential integrals\label{Sec.potentials}}

We want to extend the method of boundary layer potential to manifolds
with cylindrical ends. We begin by introducing the class of manifolds
with boundary that we plan to study in this paper.

\subsection{Submanifolds with cylindrical ends}
Let $N\subset M$ be a submanifold with boundary of a manifold with
cylindrical ends. We want to generalize the method of layer potentials
to this non-compact case.
We notice that $N$ plays a role in the method of boundary layer
potentials mostly through its boundary $\pa N$. (We shall make our assumptions
on $N$ more precise below in Definition \ref{def.m.b.c.e}.)
Because of this, we shall formulate
some of our results in the slightly more general setting when $\pa N$ is
replaced by a suitable submanifold of codimension one.

\begin{definition}\label{def.submanif}\ Let $M = M_1 \cup (\pa M_1 \times
(-\infty, 0])$ be a manifold with cylindrical ends. A submanifold with
cylindrical ends of $M$ is a submanifold $Z \subset M$ such that
$$
    Z \cap (\pa M_1 \times (-\infty, 0]) = Z' \times (-\infty, 0],
$$
for some submanifold $Z' \subset \pa M_1$. We shall write then
$Z \sim Z' \times (-\infty, 0]$.
\end{definition}

We shall fix $Z,Z'$ as above in what follows. Our main interest is
of course when $Z = \pa N$, but for certain reasonings, it is useful to
allow this slightly greater level of generality.

Let us recall from \cite[vol. II, Proposition 2.8]{TaylorB},
that a distribution $L$ on $\RR^n \times \RR^n$ is the kernel of
a classical pseudodifferential operator of order $-j$, $j = 1, 2,
\ldots, $ if, and only if,
\begin{equation}\label{eq.k.homog}
    L \sim \sum_{l = 0}^\infty (q_l(x,z) + p_l(x,z) \ln |z|)
\end{equation}
where $q_l$ are smooth functions of $x$ with values distributions in
$z$ that are homogeneous of degree $j + l - n$ and smooth for $z \not = 0$,
and $p_l$ are polynomials homogeneous of degree $j + l - n$.
(The sign ``$\sim$'' in Equation \eqref{eq.k.homog} above means that
the difference $L - \sum_{l = 0}^N (q_l(x,z) + p_l(x,z) \ln |z|)$ is
as smooth as we want if $N$ is chosen large enough.)

It is not difficult to check that the converse holds true
also for $j = 0$ under some additional conditions, for example
when $p_0 = 0$ and $q_0(x,z)$ is odd in $z$ and the associated
distribution is defined by a principal value integral.

\begin{theorem}\label{theorem.res}\
Let $M$ be a manifold with cylindrical ends and let $Z \subset M$ be a
codimension one submanifold with cylindrical ends, as in Definition
\ref{def.submanif}. If $P \in \IPS{m}$, $m < -1$, is given by the
kernel $K \in \CI(M^2 \smallsetminus M)$, then the restriction of $K$
to $Z^2 \smallsetminus Z$ extends uniquely to the kernel of an
operator $P_{Z} \in \IPSZ{m + 1}$. The same result holds true with
$\APS{m}$ and $\APSZ{m + 1}$ replacing $\IPS{m}$ and $\IPSZ{m+1}$.

Moreover, if $\sigma_m(P)$ is odd, then we can also allow $m = -1$,
provided that we define $P_{Z}$ by using a principal value integral.
\end{theorem}

\begin{proof}\
Let $P \in \IPS{m}$. Then $K$ is supported in a set of the form
$$
    V_\epsilon := \{(x,y) \in M^2, \, \text{dist}(x,y) < \epsilon\},
$$
by Definition \ref{def.t.i}.  Clearly the restriction of $K$ to $Z^2
\smallsetminus Z$ will be supported in $V_\epsilon \cap
Z^2$. Moreover, by standard (local) arguments, namely Equation
\eqref{eq.k.homog} above, $K \vert_{Z \times Z}$ is the kernel of a unique
pseudodifferential operator on $Z$ of order
$\le m + 1$.  (See \cite{TaylorA, TaylorB}). The translation
invariance of this operator follows from the definition.

To prove the same result for operators that are almost translation
invariant in a neighborhood of infinity, it is enough to do this for
order $-\infty$ operators. More precisely, we need to check
that if $T \in \APS{-\infty}$, then $T_Z \in \APSZ{-\infty}$.
This statement is local in a neighborhood
of $Z$ in the following sense. Let $\phi$ be a smooth function on $M$
that is translation invariant in a neighborhood of infinity, $\phi =
1$ in a neighborhood of $Z$ and with support in a small neighborhood
of $Z$. The statement for $T$ is equivalent to the corresponding
statement for $\phi T \phi$. We can assume then that $M = Z \times S^1$,
with $Z$ identified with $Z \times \{1\}$. By the Theorem
\ref{theorem.prod}, we can write $T = T(\theta, \theta')$, $\theta,
\theta' \in S^1$ to be a smooth function with values in
$\APSZ{-\infty}$. The result then follows because $T_Z = T(1,1)$.
\end{proof}

We need now to investigate the relation between restriction
to the submanifold $Z$ of codimension one in $M$ and indicial operators.

\begin{proposition}\label{prop.comp}\ Let $Z\subset M$ be as in Definition
\ref{def.submanif}, with $Z$ of codimension one, $Z \sim Z' \times
(-\infty,0]$, in a neighborhood of infinity. Let $P \in \APS{m}$, $m
\le -1$. Then $\tilde P_{Z' \times \RR} = \widetilde{P_{Z}}$ and
$[\hat P(\tau)]_{Z'} = \widehat{P_{Z}}(\tau)$.
\end{proposition}

\begin{proof}\ This follows from definitions, as follows. First we notice that
both statements of the Proposition are
local in a neighborhood of infinity, so we can assume that $Z = Z' \times \RR$.
The first relation then is automatic. For the second relation we also
use the fact that the restriction to $Z'$ and the Fourier transform in the
$\RR$-direction commute.
\end{proof}

\subsection{Boundary layer potential integrals}
We now proceed to define the boundary layer potential integrals. Let $M$
be a manifold with cylindrical ends and $Z \subset M$ be a submanifold
with cylindrical ends {\em of codimension one}. (Later on we shall
restrict ourselves to the case when $Z = \pa N$, where $N \subset M$
is a submanifold with boundary and cylindrical ends. For now though,
it is more convenient to continue to consider this more general case.)

Let $\delta_{Z}$ be the surface measure on $Z$, regarded as a
distribution on $M$. If $f \in L^2(Z)$, then

\begin{equation}
    f \otimes \delta_{Z} \in H^{-a}(M)\,, \quad a > 1/2.
\end{equation}

\noindent Similarly, if $\delta'_{Z}$ is the normal derivative of
$\delta_{Z}$, then

\begin{equation}
    f \otimes \delta'_{Z} \in H^{-a-1}(M)\,, \quad a > 1/2.
\end{equation}

\begin{definition}\ Fix a smooth function $V \ge 0$, $V \in \IPS{0}$,
$V$ not identical equal to $0$ on $M$. As before, we shall continue to
denote by $\Delta = \Delta_M$ the Laplace operator on $M$.  Let $f \in
L^2(Z)$ and $a > 1/2$.  The {\em single layer potential integral
associated to $Z \subset M$ and $\Delta + V$} is defined as
$$
    \Sl(f) := (\Delta + V)^{-1}(f \otimes \delta_{Z}) \in
    H^{2-a}(M),
$$
and the {\em double layer potential integral associated to $Z \subset
M$ and $V$} is defined as
$$
    \Dl(f) := (\Delta + V)^{-1}(f \otimes \delta'_{Z}) \in
    H^{1-a}(M).
$$
\end{definition}

Assume that the normal bundle of $Z$ in $M$ is oriented (so there will
be a {\em positive side} and {\em negative side} of $Z$ in $M$).  As
in \cite{MTJFA} we shall denote by $f_{\pm}$ the non-tangential limits
of some function defined on $M \smallsetminus Z$, when we approach $Z$
from the positive side ($+$), respectively from the negative side
($-$), provided, of course, that these limits exist pointwise almost
everywhere. (It is here where we need the normal bundle to $Z$ to be
oriented.)

We now begin to follow the strategy of \cite{MTJFA}.  Let

\begin{equation}\label{eq.def.S}
    S : = [(\Delta + V)^{-1}]_{Z} \in \APSZ{-1}.
\end{equation}

\noindent We shall fix in what follows a vector field $\pa_\nu$ on $M$
that is normal to $Z$ at every point of $Z$.  The principal symbol of
the order $-1$ operator $(\Delta + V)^{-1}\pa_\nu^*$ is odd, so we can
also define

\begin{equation}\label{eq.def.K}
    K : = [(\Delta + V)^{-1} \pa_\nu^*]_{Z} \in \APSZ{0}.
\end{equation}

\begin{proposition}\label{prop.sk}\ With the above notation, the operator
$S$ of Equation~\eqref{eq.def.S} is elliptic. Moreover, the zero principal
symbol of $K$ vanishes, $\sigma_0(K) = 0$, and hence actually
%\begin{equation}\label{class.K}
    $K\in \APSZ{-1}.$
%\end{equation}
\end{proposition}

\begin{proof}\ First, the fact that $S$ is elliptic follows from a symbol
calculation (which is local in nature) analogous to \cite[(3.42), p.\,33]{MMT}.
In fact, similar considerations show that $\sigma_0(K) = 0$ so, in fact,
$K\in\APSZ{-1}$. See also the discussion in
\cite[vol. II, Proposition~11.2, p.\,36]{TaylorB}.
\end{proof}

\begin{theorem}\label{theorem.map.slp}\ Let $Z \subset M$ be a codimension
one submanifold with cylindrical ends. Assume the normal bundle to $Z$ is
oriented. Given $f \in L^2(Z)$, we have
$$
    \Sl(f)_+ = \Sl(f)_- = S f
$$
as pointwise a.e. limits. Also, using the notation of Equation
\eqref{eq.def.K} above, we have
$$
    \pa_\nu \Sl(f)_{\pm} = \Bigl(\pm\ha I + K^*\Bigr)f,
$$
where $K^*$ is the formal transpose of $K$.
\end{theorem}

\begin{proof}\ Let us write $T:= (\Delta + V)^{-1} = P + R$, where
$P\in \IPS{m}$  (so it is translation invariant in a
neighborhood of infinity) and $R \in \APS{-\infty}$.  The first
statement of the proposition, namely
$$
    [T(f\otimes \delta_Z)]_{\pm } = T_Z f
$$
is clearly linear in $T \in \APS{m}$, $m < -1$. It is enough
then to prove it for $P$ and $R$ separately.

For $T = (\Delta + V)^{-1}$ replaced by $P$, this is a local statement
(because $P$ is properly supported), which then follows from
\cite[Proposition 3.8]{MTJFA}.

For $T$ replaced by $R$, we argue as in the proof of
Theorem \ref{theorem.res} that we can assume that $M = Z \times S^1$,
with $Z$ identified with the submanifold $Z \times \{1\}$.
Then we use again Theorem \ref{theorem.prod} to write $R = R(\theta,\theta')$,
for some smooth function with values in $\APSZ{-\infty}$.

This gives
$$
    R(f \otimes \delta_Z) (z,\theta) = [R(\theta,1)f](z)
$$
and $R_Z = R(1,1)$. Let $g_\theta(z) = R(f \otimes \delta_Z)(z,\theta)$.
The assumptions on the function $R(\theta,\theta')$
% $T(\theta,\theta')$
guarantee that the function
$$
    S^1\ni\theta\mapsto g_\theta \in H^m(M)
$$
is continuous (in fact, even $\CI$) for any $m$. Then
$$
    [R(f\otimes \delta_Z)]_{\pm}=\displaystyle{\lim_{\theta\to
    1\pm 0}}g_\theta =g_1:= R(1,1)f = R_Zf.
$$
\end{proof}

The following theorem is proved in a completely similar way, following
the results of \cite[Proposition 3.8]{MTJFA}.

\begin{theorem}\label{theorem.map.dlp}\ Let $Z$ be a {\em codimension one
submanifold with cylindrical ends of $M$ with oriented normal bundle}.
Given $f \in L^2(Z)$, we have
$$
    \Dl(f)_{\pm} =\Bigl(\mp\ha I + K\Bigr)f
$$
as pointwise a.e. limits.
\end{theorem}

We can replace the pointwise almost everywhere limits with
$L^2$--limits both for the tangential limits of the single and
double layer potentials; see Theorem \ref{theorem.h.reg2}.

For further reference, let us discuss now the ``trace theorem'' for
codimension one submanifolds in our setting. See \cite{ALNV} for more
details and results of this kind for manifolds with a Lie structure at
infinity.

\begin{proposition}\ Let $Z \subset M$ be a submanifolds with cylindrical ends
of the manifold with cylindrical ends $M$. Then the restriction map
$\CIc(M) \to \CIc(Z)$ extends to a continuous map $H^s(M) \to
H^{s-1/2}(Z)$, for any $s > 1/2$.
\end{proposition}

\begin{proof}\
We can assume, as in the proof of Theorem \ref{theorem.map.slp}, that
$M = Z \times S^1$. Since the Sobolev spaces $H^s(M)$ and
$H^{s-1/2}(Z)$ do not depend on the metric on $M$ and $Z$, as long as
these metrics are compatible with the structure of manifolds with
cylindrical ends, we can assume that the circle $S^1$ is given the
invariant metric making it of length $2\pi$ and that $M$ is given the
product metric.

Then $\Delta = \Delta_Z + \Delta_{S^1}$ and $\Delta_{S^1}=-
\pa_\theta^2$ has spectrum $\{ 4\pi^2n^2 \}$, $n \in \ZZ$. We can
decompose $L^2(Z\times S^1)$ according to the eigenvalues $n \in \ZZ$
of $(2\pi \imath)^{-1} \pa_{\theta}$:
% $\frac{\imath\pa}{2\pi\pa\theta}$:
$$
    L^2(Z \times S^1) \simeq \oplus_{n \in \ZZ} L^2(Z \times
    S^1)_n \simeq \oplus_{n \in \ZZ} L^2(Z),
$$
where the isomorphism $L^2(Z\times S^1)_n\simeq\oplus_{n \in \ZZ}
L^2(Z)$ is obtained by restricting to $1 \in S^1$.

To prove our proposition, it is enough to check that if $\xi_n \in
L^2(Z)$ is a sequence such that
\begin{equation}\label{eq.sum.n}
    \sum_n \|( 1 + n^2 + \Delta_Z)^{s/2}\xi_n\|^2 < \infty
\end{equation}
then $\sum (1 + \Delta_Z)^{s/2-1/4}\xi_n$ is convergent.

Let $C = 1 + \int_\RR (1 + t^2)^{-s}dt$ and assume that each $\xi_n$
is in the spectral subspace of $\Delta_Z$ corresponding to
$[m,m+1)\subset\RR_+$.  Then
$$
    (1 + m^2)^{s - 1/2} \Bigl( \sum_n \| \xi_n \| \Bigr )^2 \le C
    \sum_n \| (1 + n^2 + m^2)^{s/2} \xi_n\|^2.
$$
Since the constant $C$ is independent of $m$ and the spectral spaces
of $\Delta_Z$ corresponding to $[m,m+1) \subset \RR$ give an
orthogonal direct sum decomposition of $L^2(Z)$, this checks Equation
\eqref{eq.sum.n} and completes the proof.
\end{proof}

\subsection{Higher regularity of the layer potentials}
We shall not need the following results in what follows. We include
them for completeness and because they give a better intuitive picture
of the properties of layer potentials.  Choose a small open tubular
neighborhood $U$ of $Z$ in $M$, such that $U \simeq Z \times
(-\epsilon, \epsilon)$ via a diffeomorphism that is compatible with
the cylindrical ends structure of $Z$ and $M$. For example, assume
that $\pa_\nu$ is a vector field on $M$ that is normal to $Z$ and
translation invariant in a neighborhood of infinity. Denote by
$\exp(t\pa_\nu)$ the one-parameter group of diffeomorphisms generated
by $\pa_\nu$. (This group exists because $\pa_{\nu}$ extends to the
canonical compactification of $M$ to a manifold with boundary $\simeq
M_1$.)  Then the range $U = U_\epsilon$ of the map
$$
    Z \times (-\epsilon,\epsilon) \ni (z,t) \mapsto \Psi(z,t) :=
    \exp(t\pa_\nu)z \in M
$$
is a good choice, for $\epsilon > 0$ small enough. In particular, for
$\epsilon$ small enough, the complement $U_\epsilon^c$ of $U_\epsilon$
is a smooth submanifold with boundary, such that its boundary $\pa
U_\epsilon^c$ is a submanifold with cylindrical ends. Moreover, $\pa
U_\epsilon^c = Z_{-\epsilon} \cup Z_{+\epsilon}$ is the disjoint union
of two manifolds diffeomorphic to $Z$ via $Z \simeq Z \times \{ \pm
\epsilon \} \simeq Z_{\pm \epsilon}$, where the second map is given by
$\Psi$.

Denote by $H^m(U_\epsilon^c)$ the space of restrictions of
distributions in $H^m(M)$ to (the interior of) the complement of
$U_\epsilon$.

The following two theorems describe the mapping properties of the
single and double layer potentials. Since the statements and proofs
work actually in greater generality, we begin with some more general
results, which we shall then specialize to the case of single and
double layer potentials.

\begin{theorem}\label{theorem.int.reg}\
Let $U \simeq Z \times (-\epsilon, \epsilon)$ be a tubular
neighborhood of $Z$ in $M$ (as above) and let $T \in
\APS{m}$. Restriction to $U^c$ defines for any $s$ continuous maps
$$
    H^s(Z) \ni f \mapsto T( f \otimes \delta_Z) \in
    H^{\infty}(U^c),
$$
which are translation invariant in a neighborhood of infinity, for any
tubular neighborhood $U$ of $Z$.
\end{theorem}

\begin{proof}\ Let $\psi_0$ and $\psi_1$ be smooth functions on $M$
and $T \in \APS{m}$. Assume the following:\ \ $\psi_0$ and $\psi_1$
are translation invariant in a neighborhood of infinity;\ $\psi_0$ is
equal to $1$ in a neighborhood of $Z$;\ $\psi_1$ vanishes in a
neighborhood of the support of $\psi_0$; and $\psi_0$ is equal to $1$
in a neighborhood of $U^c$. Then
$$
    T(f\otimes\delta_Z) \vert_{U^c} = (\psi_1 T \psi_0) (f\otimes
     \delta_Z)
$$
and $\psi_1 T \psi_0 \in \APS{-\infty}$ because the supports of
$\psi_0$ and $\psi_1$ are disjoint.
\end{proof}

Consider now $U = U_\epsilon \simeq Z \times (-\epsilon, \epsilon)$,
for $\epsilon > 0$ small enough, where the last diffeomorphism is
given by the exponential map. Then decompose $\pa U_\epsilon^c =
Z_{+\epsilon} \cup Z_{-\epsilon}$ as a disjoint union, as above. In
particular, we fix the diffeomorphisms $Z \simeq Z_{\pm \epsilon}$
defined by the exponential, as above. Then the traces of the
restrictions to $U_{\epsilon}^{c}$
\begin{equation}\label{eq.def.Tt}
    H^s(Z) \ni f \to T_{\pm \epsilon}f :=
    T(f \otimes \delta_Z) \vert_{Z_{\pm \epsilon}} \in
    H^{s'}(Z_{\pm \epsilon}) \simeq H^{s'}(Z)
\end{equation}
define continuous operators $T_{\pm \epsilon} : H^s(Z) \to H^{s'}(Z)$, for any
$s, s' \in \RR$.

We fix in what follows $\epsilon > 0$ as above. Similarly, we obtain
operators $T_{\pm t} : H^{s}(Z) \to H^{s'}(Z)$, for any $t \in
(0,\epsilon]$ and any $s, s' \in \RR$.

\begin{theorem}\label{theorem.h.reg}\ Let $T \in \APS{m}$ and $T_{t}$ be
as above, Equation \eqref{eq.def.Tt}. Then $T_{\pm t} \in
\APSZ{-\infty}$ and the two functions
$$
    (0,\epsilon] \ni t \to t^l \pa_t^k T_{\pm t} \in \APSZ{m + 1 +
    k - l + \delta}
$$
extend by continuity to $[0,\epsilon]$ if $\delta > 0$.
These extensions are bounded for $\delta = 0$.
\end{theorem}

\begin{proof}\ The proof is based on the ideas in
\cite[vol. II, Ch. 7, Sec. 12]{TaylorB}, especially Theorem 12.6,
and some local calculations. Here are some details.

Since the statement of the theorem is ``linear'' in $T$, it is enough
to prove it for $T \in \PS{m}$ and for $T \in \APS{-\infty}$. The
later case is obvious -- in fact, it is already contained in the proof
of Theorem \ref{theorem.map.slp}. Then, we can further reduce the
proof to the case when $T = s_0(T_1)$, with $T_1 \in \Psi_{\rm ai}(\pa
M_1 \times \RR)^{\RR}$, and to the case when $T$ has compactly
supported Schwartz kernel.  Again, the second case is easier, being an
immediate consequence of the corresponding result for the compact case.
Because the second case involves a similar argument, we shall
nevertheless discuss this here.

Assume, for the next argument, that $M$ is compact. Since the result
is true for regularizing operators, we can use a partition of unity to
localize to the domain of a coordinate chart. This allows then to
further replace $M$ with $\RR^n$, $Z$ with $\RR^{n-1}$, and $T$ with
an operator of the form $T = a(x,D)$, with $a(\;,\;)$ in H\"ormander's
symbol class $S^{m}_{1,0} = S^{m}_{1,0}(\RR^n)$
\cite[vol. II]{TaylorB} of functions that satisfy uniform estimates in
the space variable $x$ (and the usual symbolic estimates in the dual
variable).

Let $(x', x_n) \in \RR^{n-1} \times \RR$ and $(\xi', \xi_n) \in
\RR^{{n-1}\ast}\times \RR^*$ be the usual decomposition of the variables.
Also, let
$$
    a_{t}(x',\xi') = (2\pi \imath)^{-1} \int_{\RR} e^{\imath t
    \xi_n} a(x, t, \xi', \xi_n) d\xi_n.
$$
Then $a_{t}$ is such that $T_{t} = a_{t}(x,D)$ and the (two) functions
$t^l \pa_{t}^k a_{\pm t}$ extend to continuous functions $[0,\epsilon]
\to S^{m + 1 + l - k + \delta}_{1,0}$, for any $\delta > 0$.  These
extensions are bounded as functions with values in $S^{m + 1 + l -k}_{1,0}$.
This completes the proof of our result for the case $M$ compact.

Let us consider now to the case when $T = s_0(T_1)$. We can assume
that $M = \pa M_1 \times \RR$ and that $T$ is $\RR$--invariant.  The
proof is then the same as in the case $M$ compact, but using local
coordinates on $\pa M_1$ instead of on $M$, and making sure that all
our symbols and all maps preserve the $\RR$-invariance. This completes
the proof of our result.
\end{proof}

A consequence of the above theorem is the following continuity result.

\begin{corollary}\label{cor.cont}\ Let $T \in \APS{-1}$.
\begin{enumerate}[(i)]
\item\ If $f \in H^{m}(Z)$, then the functions $(0,\epsilon] \ni t
\mapsto T_{\pm t} f \in H^{m}(Z)$ extend by continuity at $0$.
\item\ If $f \in H^{\infty}(Z)$, then the mappings $(0, \epsilon]
\times Z \mapsto \big( T_{\pm t}f \big) (z)$ extend to functions in
$H^{\infty}([0,\epsilon] \times Z)$.
\end{enumerate}
\end{corollary}

\begin{proof}\
Denote by $\mathcal L(X,Y)$ the normed space of bounded operators
between two Banach spaces $X$ and $Y$. Theorem \ref{theorem.h.reg}
ensures that $(0, \epsilon] \to T_{\pm t}$ have continuous extensions
to functions
$$
    [0,\epsilon] \to \mathcal L(H^{m + \delta}(Z), L^2(Z)),
$$
for $\delta > 0$. For $\delta = 0$ these extensions will be bounded.

This proves the first part of our result as follows. If $f \in H^{m +
\delta}(Z)$, then the functions $T_{\pm t}f \in L^2(Z)$ extend by
continuity on $[0, \epsilon]$ because $T_{\pm t}$ extend by continuity
on $[0, \epsilon]$ as maps to $\mathcal L(H^{m + \delta}(Z), L^2(Z))$.
Since $H^{m + \delta}(Z)$, $\delta > 0$, is dense in $H^{m}(Z)$ and
$T_{\pm t}$ are bounded as maps $[0,\epsilon] \to \mathcal L(H^{m}(Z),
L^2(Z))$, the result follows from an $\epsilon/3$--type argument.

To prove (ii), it is enough to prove then that $\pa_{t}^b (I +
\Delta_{Z})^a T_{\pm t}f$ is in $L^{2}([0, \epsilon] \times Z)$, for
any $a , b \in \NN$. Using again Theorem \ref{theorem.h.reg}, we know
that $\pa_{t}^b (I + \Delta_{Z})^a T_{\pm t}$ extend to continuous
functions $[0, \epsilon] \to \APS{c}$, with $c = m + 2 + k + a$, (take
$\delta = 1$).  Since $f \in H^{\infty}(Z) \subset H^{c}(Z)$, the
functions $(0, \epsilon] \ni t \to \pa_{t}^b (I + \Delta_{Z})^a T_{\pm
t}f \in L^2(Z)$ extend by continuity to a function defined on
$[0,\epsilon]$. This extension is then in $L^2([0,\epsilon] \times
Z)$.
\end{proof}

We can specialize all the above results to $T = (\Delta + V)^{-1}$ or
$T = (\Delta + V)^{-1}\pa^*_\nu$. This gives maps $S_{\pm t}(f) :=
\Sl(f)\vert_{Z_{\pm t}}$ and $D_{\pm t}(f) := \Dl(f) \vert_{Z_{\pm t}}$,
where $t \in (0, \epsilon]$.

\begin{theorem}\label{theorem.h.reg2}\ Using the notation we have
just introduced, we have
\begin{enumerate}[(i)]
\item\ $S_{\pm t}, D_{\pm t} \in \APSZ{-\infty}$ and the
functions $(0,\epsilon] \ni t \to t^l\pa_t^k S_{\pm t} \in
\APSZ{\delta -1 + k - l}$ and $(0,\epsilon] \ni t \to
t^l \pa_t^k D_{\pm t} \in \APSZ{\delta + k - l}$
extend by continuity to $[0,\epsilon]$, for $\delta > 0$. For $\delta = 0$
these functions are bounded.
\item\ If $f \in L^2(Z)$, then the functions $t \to S_{\pm t}f, \;
D_{\pm t}f \in L^2(Z)$ extend by continuity to $[0,\epsilon]$.
\item\ If $f \in H^{\infty}(Z)$, then the restrictions of $\Sl(f)$
and $\Dl(f)$ to $Z \times [-\epsilon, 0)$ and, respectively,
$Z\times (0,\epsilon]$ extend to functions in
$H^{\infty}(Z \times [-\epsilon,0])$,
respectively in $H^{\infty}(Z \times [0, \epsilon])$.
\end{enumerate}
\end{theorem}

\section{Layer potentials depending on a parameter\label{Sec.parameter}}

The aim of this section is to investigate the invertibility
of layer potential operators which depend on a parameter
$\tau\in\RR$, via a method initially developed by G. Verchota
in \cite{Ve}, for the case of the flat-space Laplacian.
The novelty here is to derive estimates which
are {\it uniform} with respect to the real parameter $\tau$.

Let $\manif$ be a smooth, {\it compact, boundaryless} Riemannian
manifold, and fix a reasonably regular subdomain $\Omega\subset\manif$
(Lipschitz will do). Here, $\manif$ will play the role of $\pa M_1$ in
our standard notation and, anticipating notation introduced in the next
section, $\Omega$ will play the role of the {\em exterior} of $X$.

Set $\nu$ for the outward unit conormal to
$\Omega$ and $d\sigma$ for the surface measure on $\partial\Omega$
(naturally inherited from the metric on $\manif$). The departure point
is the following Rellich type identity:

\begin{multline} \label{eq.3.3.1}
    \int_{\partial\Omega}\langle\nu,w\rangle
    \bigl\{|\nabla_{tan}u|^2-|\partial_\nu u|^2\bigr\}\,d\sigma \\
    \qquad\quad =2\mbox{Re}\,\int\limits_{\partial\Omega} \langle
    w_{tan},\nabla u\rangle\partial_\nu\bar{u}\,d\sigma
    -2\mbox{Re}\,\int_{\Omega}\langle\nabla\bar{u},w\rangle
        \Delta_\manif\,u\,dx \\
        \qquad\qquad
        +\mbox{Re}\,\int_{\Omega}\bigl\{({\rm div}\,w)|\nabla u|^2
        -2(\mathcal{L}_wg)(\nabla u,\nabla\bar{u})\bigr\}\,dx,
\end{multline}

\noindent which, so we claim, is valid for a (possibly complex-valued) scalar
function $u$ and a real-valued vector field $w$ (both sufficiently
smooth, otherwise arbitrary) in $\Omega$.  Hereafter, the subscript
`tan' denotes the tangential component relative to $\partial\Omega$.
At the level of vector fields, $\nabla$ is used to denote the
Levi-Civita connection on $\manif$.
Also, $\mathcal{L}_w g$ stands for the Lie derivative of the metric
tensor $g$ with respect to the field $w$; recall that, in general,
$$
\mathcal{L}_w g(X,Y)=\langle\nabla_Xw,Y\rangle+\langle\nabla_Yw,X\rangle,
$$
for any two vector fields $X,Y$.

To prove (\ref{eq.3.3.1}), consider the vector field
$F:=|\nabla u|^2w-2(\partial_w u)\nabla\bar{u}$ and compute

\begin{eqnarray}\label{new1}
\langle\nu,F\rangle
& = & |\nabla u|^2\langle\nu,w\rangle
-2\,(\partial_w u)(\partial_{\nu}\bar{u})
\nonumber\\[4pt]
& = & |\nabla u|^2\langle\nu,w\rangle
-2\,\langle w_{\rm tan},\nabla u\rangle\,\partial_{\nu}\bar{u}
-2\,|\partial_\nu u|^2\langle\nu,w\rangle
\nonumber\\[4pt]
& =& \langle\nu,w\rangle
\Bigl(|\nabla_{\rm tan} u|^2-|\partial_\nu u|^2\Bigr)
-2\langle w_{\rm tan},\nabla u\rangle\,\partial_{\nu}\bar{u},
\end{eqnarray}

\noindent by decomposing $w=w_{\rm tan}+\langle\nu,w\rangle\nu$
and $|\nabla u|^2=|\nabla_{\rm tan} u|^2 +|\partial_\nu u|^2$.
Furthermore,

\begin{equation}\label{new3}
{\rm div}\,F = ({\rm div}\,w)|\nabla u|^2+w(|\nabla u|^2)
-2\,(\partial_w u)\Delta_{\manif}\,\bar{u}-2\,\nabla u(\partial_w\bar{u}).
\end{equation}

\noindent Given the current goal, the first and the third terms suit
our purposes; for the rest we write

\begin{multline}%\label{new2}
    w(|\nabla u|^2)-2\,\nabla u(\partial_w\bar{u})  =  w(|\nabla
    u|^2)-2\,\nabla u(w(\bar{u})) \nonumber \\[4pt]  =
    w(|\nabla u|^2)-2\,[\nabla u,w]\bar{u}-2\,w(\nabla u(\bar{u}))
    \nonumber \\[4pt]  =  w(|\nabla
    u|^2)+2\langle\nabla_w(\nabla u),\nabla\bar{u}\rangle
    -2\langle\nabla_{\nabla u}w,\nabla\bar{u}\rangle -2\,w|\nabla
    u|^2 \nonumber \\[4pt]  =  -w(|\nabla u|^2)+{\rm
    Re}\,[w(|\nabla u|^2)] -2({\mathcal{L}}_w\,g)(\nabla
    u,\nabla\bar{u}),
\end{multline}

\noindent where the third equality utilizes the fact that $\nabla$ is
torsion-free. Since the real parts of the first two terms in the last
expression above cancel out, it ultimately follows that

\begin{equation}\label{new4}
    {\rm Re}\,({\rm div}\,F)= ({\rm div}\,w)|\nabla u|^2
    -2\,{\rm Re}\,[(\partial_w\bar{u})\Delta_{\manif}\,u]
    -2\,{\rm Re}\,({\mathcal{L}}_w\,g)(\nabla u,\nabla\bar{u}).
\end{equation}

\noindent Thus, the Rellich identity (\ref{eq.3.3.1}) follows from
\eqref{new4}, \eqref{new1}, and the Divergence Theorem, after
taking the real parts.

Another general identity (in fact, a simple consequence of the
Divergence Theorem) that is useful here is

\begin{equation}\label{eq.3.3.2}
    \int_{\partial\Omega}|u|^2\langle w,\nu\rangle\,d\sigma
    =\mbox{Re}\,\int_\Omega\{2u\langle\nabla\bar{u},w\rangle
    +({\rm div}\,w)|u|^2\}\,dx.
\end{equation}

To proceed, fix a nonnegative scalar potential $W\in C^\infty(\manif)$
and for the remainder of this subsection assume that

\begin{equation}\label{eq.3.3.3}
    (\Delta_\manif + \tau^2 + W)u=0 \quad \mbox{in} \quad \Omega,
\end{equation}

\noindent where $\tau\in\mathbb{R}$ is an arbitrary parameter (fixed for
the moment). Our immediate objective is to show that

\begin{equation}\label{eq.3.3.4}
    \int_{\partial\Omega}|\partial_\nu u|^2\,d\sigma \leq
    C\int_{\partial\Omega}\{|\nabla_{tan}u|^2+(1+\tau^2)|u|^2\}\,d\sigma,
\end{equation}

\noindent uniformly in $\tau$, and that for each $\varepsilon>0$ there
exists a finite constant $C=C(\Omega,\varepsilon)>0$ so that

\begin{equation}\label{eq.3.3.5}
    \int_{\partial\Omega}\{|\nabla_{tan}u|^2+\tau^2|u|^2\}\,d\sigma
    \leq C\int_{\partial\Omega}|\partial_\nu u|^2\,d\sigma
    +\varepsilon\int_{\partial\Omega}|u|^2\,d\sigma
\end{equation}

\noindent uniformly in the parameter $\tau\in\mathbb{R}$. We shall also
need a strengthened version of (\ref{eq.3.3.5}) to the effect that

\begin{equation}\label{eq.3.3.6}
    W>0 \mbox{ in } \Omega \Longrightarrow
    \int_{\partial\Omega}\{|\nabla_{tan}u|^2+(1+\tau^2)|u|^2\}\,d\sigma
    \leq C\int_{\partial\Omega}|\partial_\nu u|^2\,d\sigma
\end{equation}

\noindent uniformly in the parameter $\tau\in\mathbb{R}$.

With an eye on (\ref{eq.3.3.5}), let us recall Green's first identity
for the function $u$ that we assumed to satisfy Equation \eqref{eq.3.3.3}

\begin{equation*}%\label{eq.3.3.7}
\int_\Omega\{|\nabla u|^2+\tau^2|u|^2+W|u|^2\}\,dx=
\mbox{Re}\,\int_{\partial\Omega}\bar{u}\,\partial_\nu u\,d\sigma
\end{equation*}

\noindent which readily yields the energy estimate

\begin{equation}\label{eq.3.3.8}
\int_\Omega\{|\nabla u|^2+\tau^2|u|^2+W|u|^2\}\,dx
\leq \int_{\partial\Omega}|u||\partial_\nu u|\,d\sigma.
\end{equation}

\noindent In turn, this further entails

\begin{equation}\label{eq.3.3.9}
\int_\Omega\tau^2|\nabla u||u|\,dx\leq C|\tau|\,
\int_\Omega\{\tau^2|u|^2+|\nabla u|^2\}\,dx\leq C|\tau|\,
\int_{\partial\Omega}|u||\partial_\nu u|\,d\sigma,
\end{equation}

\noindent uniformly in $\tau$.

Let us now select $w$ to be transversal to $\partial\Omega$, i.e.

\begin{equation}\label{eq.3.3.10}
{\rm ess\,inf}\,\langle w,\nu\rangle>0\mbox{ on }\partial\Omega,
\end{equation}

\noindent something which can always be arranged given that
$\partial\Omega$ is assumed to be Lipschitz. This, in concert
with (\ref{eq.3.3.2}), then gives

\begin{equation}\label{eq.3.3.11}
\int_{\partial\Omega}|u|^2\,d\sigma\leq
C\int_{\Omega}\{|u|^2+|\nabla u||u|\}\,dx.
\end{equation}

\noindent Multiplying (\ref{eq.3.3.11}) with $\tau^2$ and then invoking
(\ref{eq.3.3.8})-(\ref{eq.3.3.9}) eventually justifies the estimate

\begin{equation}\label{eq.3.3.12}
\int_{\partial\Omega}\tau^2|u|^2\,d\sigma\leq
C\int_{\partial\Omega}(1+|\tau|)|\partial_\nu u||u|\,d\sigma.
\end{equation}

\noindent Next, make the (elementary) observation that
for every $\varepsilon,\delta>0$ there exists $C=C(\varepsilon,\delta)>0$
so that

\begin{equation}\label{eq.3.3.13}
    (1+|\tau|)|\partial_\nu u||u|\leq \delta\tau^2|u|^2
    +C|\partial_\nu u|^2+\varepsilon|u|^2,
\end{equation}

\noindent uniformly in $\tau$. When considered in the context of
(\ref{eq.3.3.12}), the boundary integral produced by the first term in the
right side of (\ref{eq.3.3.13}) can be absorbed in the left side of
(\ref{eq.3.3.12}), provided $\delta$ is sufficiently small. Thus, with
this alteration in mind, (\ref{eq.3.3.12}) becomes

\begin{equation*}\label{eq.3.3.14}
    \int_{\partial\Omega}\tau^2|u|^2\,d\sigma\leq
    C\int_{\partial\Omega}|\partial_\nu u|^2\,d\sigma
    +\varepsilon\int_{\partial\Omega}|u|^2\,d\sigma
\end{equation*}

\noindent which is certainly in the spirit of \eqref{eq.3.3.5}. In fact,
in order to fully prove the latter estimate, there remains to control the
tangential gradient in a similar fashion. To this end, observe that
\eqref{eq.3.3.10} and Rellich's identity \eqref{eq.3.3.1} give

\begin{equation*}\label{eq.3.3.15}
    \int_{\partial\Omega}|\nabla_{tan}u|^2\,d\sigma\leq
    C\int_{\partial\Omega}|\partial_\nu u|^2\,d\sigma
    +C\int_{\Omega}\tau^2|\nabla u||u|\,dx
    +C\int_{\Omega}|\nabla u|^2\,dx,
\end{equation*}

\noindent uniformly in $\tau$. With this at hand, the same type of estimates
employed before can be used once again to further bound the solid integrals
in terms of (suitable) boundary integrals. The bottom line is that

\begin{equation}\label{eq.3.3.16}
    \int_{\partial\Omega}|\nabla_{tan}u|^2\,d\sigma\leq
    C\int_{\partial\Omega}|\partial_\nu u|^2\,d\sigma
    +\varepsilon\int_{\partial\Omega}|u|^2\,d\sigma,
\end{equation}

\noindent uniformly in $\tau$, and (\ref{eq.3.3.5}) follows.

It is now easy to prove \eqref{eq.3.3.6}, having disposed off
(\ref{eq.3.3.5}). One useful ingredient in this regard is

\begin{equation}\label{eq.3.3.17}
    \int_\Omega|u|^2\,dx\leq C\int_\Omega\{|\nabla u|^2+W|u|^2\}\,dx,
\end{equation}

\noindent itself a version of Poincar\'e's inequality.
When used in conjunction with (\ref{eq.3.3.8}) and
(\ref{eq.3.3.11}), this readily yields

\begin{equation}\label{eq.3.3.18}
    \int_{\partial\Omega}|u|^2\,d\sigma
    \leq C\int_\Omega\{|\nabla u|^2+W|u|^2\}\,dx
    \leq C\int_{\partial\Omega}|u||\partial_\nu u|\,d\sigma
\end{equation}

\noindent so that, ultimately,

\begin{equation}\label{eq.3.3.19}
\int_{\partial\Omega}|u|^2\,d\sigma
\leq C\int_{\partial\Omega}|\partial_\nu u|^2\,d\sigma,
\end{equation}

\noindent in the case we are currently considering.
In concert with \eqref{eq.3.3.5}, this concludes the proof of
\eqref{eq.3.3.6}. Let us now turn our attention to the estimate
\eqref{eq.3.3.4}. For starters, Rellich's identity \eqref{eq.3.3.1}
can also be employed, along with the condition \eqref{eq.3.3.10}, to produce

\begin{equation}\label{eq.3.3.20}
    \int_{\partial\Omega}|\partial_\nu u|^2\,d\sigma\leq
    C\int_{\partial\Omega}|\nabla_{tan} u|^2\,d\sigma
    +C\int_{\Omega}\tau^2|\nabla u||u|\,dx
    +C\int_{\Omega}|\nabla u|^2\,dx,
\end{equation}

\noindent uniformly in $\tau$. Then, much as before,

\begin{eqnarray}%\label{eq.3.3.21}
    \int_{\Omega}\{\tau^2|\nabla u||u|+|\nabla u|^2\}\,dx
    & \leq & C\int_{\partial\Omega}(1+|\tau|)|\partial_\nu u||u|\,d\sigma
    \nonumber \\[4pt]
    & \leq & \delta\int_{\partial\Omega}|\partial_\nu u|^2\,d\sigma
    +C\int_{\partial\Omega}(1+\tau^2)|u|^2\,d\sigma,
\end{eqnarray}

\noindent where $\delta>0$ is chosen small and $C$ depends only on
$\Omega$ and $\delta$. With these two estimates at hand, the endgame
in the proof of (\ref{eq.3.3.4}) is clear.

After these preliminaries, we can finally address the main theme of
this subsection. More concretely, for each $\tau\in\mathbb{R}$, let
$S_\tau$, $K_\tau$ be, respectively, the single and the double layer
potential operators associated with $\Delta_\manif+\tau^2+W$ on
$\partial\Omega$ (recall that the potential $W$ was first introduced
in connection with (\ref{eq.3.3.3})). {}From the work in \cite{MTJFA},
it is known that if $\Omega$ has a Lipschitz boundary then both

\begin{equation*}%\label{eq.3.3.22}
    S_\tau:L^2(\partial\Omega)\longrightarrow H^1(\partial\Omega)
    \quad\mbox{ and }\quad
    \ha I+K_\tau:L^2(\partial\Omega)\longrightarrow L^2(\partial\Omega)
\end{equation*}

\noindent are invertible operators for each $\tau\in\mathbb{R}$.
Our objective is to study  how the norms of their inverses depend on the
parameter $\tau$. To discuss this issue, for each $\tau\in\mathbb{R}$ and
$f\in H^1(\partial\Omega)$, set

\begin{equation}\label{eq.3.3.23}
\|f\|_{H^1_\tau(\partial\Omega)}:=\|f\|_{H^1(\partial\Omega)}+
|\tau|\|f\|_{L^2(\partial\Omega)}.
\end{equation}

\noindent Thus, $\mathbb{R}\ni\tau\mapsto
\|\cdot\|_{H^1_\tau(\partial\Omega)}$ is a one-parameter family of
equivalent norms on the Sobolev space $H^1(\partial\Omega)$. The main
result of this subsection is as follows.

\vskip 0.10in

\begin{proposition}\label{prop.3.3.1}
Assume that $\Omega$ is a fixed, Lipschitz subdomain of $\manif$, and
retain the notation introduced above.  Then there exits a finite
constant $C=C(\partial\Omega)>0$, depending exclusively on the
Lipschitz character of $\Omega$, such that for each
$\tau\in\mathbb{R}$, we have

\begin{equation}\label{eq.3.3.24}
\|S_\tau^{-1}f\|_{L^2(\partial\Omega)}\leq C\|f\|_{H^1_\tau(\partial\Omega)}
\end{equation}

\noindent uniformly for $f\in H^1(\partial\Omega)$.

Furthermore, if $W>0$ on a set of positive measure in $\Omega$, then for
any $\tau\in\mathbb{R}$ we also have

\begin{equation}\label{eq.3.3.25}
\|(\ha I+K_\tau)^{-1}f\|_{L^2(\partial\Omega)}\leq
C\|f\|_{L^2(\partial\Omega)},
\end{equation}

\noindent uniformly for $f\in L^2(\partial\Omega)$.
\end{proposition}

\begin{proof} Consider first (\ref{eq.3.3.25}). Let $\Omega_+:=\Omega$,
$\Omega_-:=\manif\setminus\bar{\Omega}$, and for $f\in L^2(\partial\Omega)$,
set $u:=\mathcal{S}f$ in $\Omega_\pm$. Thus,

\begin{equation}\label{eq.3.3.26}
(u)_{+}=(u)_{-},\quad
(\nabla_{tan}u)_{+}=(\nabla_{tan}u)_{-},\quad
(\partial_{\nu}u)_{\pm}=(\pm\ha I+K_\tau^*)f.
\end{equation}

\noindent In turn, (\ref{eq.3.3.26}), (\ref{eq.3.3.4})
and (\ref{eq.3.3.6}) allow us to write

\begin{multline*} %\label{eq.3.3.27}
    \|(-\ha I+K_\tau^*)f\|_{L^2(\partial\Omega)}
    =\|(\partial_\nu u)_{-}\|_{L^2(\partial\Omega)}
    \\  \leq  C\|(u)_{-}\|_{H^1_\tau(\partial\Omega)}
    =C\|(u)_{+}\|_{H^1_\tau(\partial\Omega)}
\leq C\|(\partial_\nu u)_{+}\|_{L^2(\partial\Omega)}
    \\ =C\|(\ha I+K_\tau^*)f\|_{L^2(\partial\Omega)}.
\end{multline*}

\noindent Consequently,

\begin{multline}\label{eq.3.3.28}
    \|f\|_{L^2(\partial\Omega)}
    \leq \|(-\ha I+K_\tau^*)f\|_{L^2(\partial\Omega)}
    +\|(\ha I+K_\tau^*)f\|_{L^2(\partial\Omega)} \\
    \leq  C\|(\ha I+K_\tau^*)f\|_{L^2(\partial\Omega)}
\end{multline}

\noindent for some constant $C=C(\partial\Omega)>0$ independent of $\tau$.
Going further, if $\mathcal{L}(X):=\mathcal{L}(X,X)$, the normed algebra
of all bounded operators on a Banach space $X$, then (\ref{eq.3.3.28}) entails

\begin{equation*}%\label{eq.3.3.29}
\|(\ha I+K_\tau)^{-1}\|_{\mathcal{L}\bigl(L^2(\partial\Omega)\bigr)}
=\|(\ha I+K_\tau^*)^{-1}\|_{\mathcal{L}\bigl(L^2(\partial\Omega)\bigr)}
\leq C.
\end{equation*}

\noindent This takes care of (\ref{eq.3.3.25}).

As for (\ref{eq.3.3.24}), the argument is rather similar, the main step
being the derivation of the estimate

\begin{equation*}%\label{eq.3.3.30}
    \|f\|_{L^2(\partial\Omega)}\leq
    C\|\nabla_{tan}(S_\tau f)\|_{L^2(\partial\Omega)}
    +C(1+|\tau|)\|S_\tau f\|_{L^2(\partial\Omega)},
\end{equation*}

\noindent out of (\ref{eq.3.3.26}) and (\ref{eq.3.3.4}), when the latter
is written both for $\Omega_+$ and $\Omega_-$. Once again, the crux of the
matter is that the intervening constant $C=C(\partial\Omega)>0$ is
independent of $\tau$. The proof is finished.
\end{proof}

\section{The Dirichlet problem\label{Sec.Dirichlet}}

We now apply the results we have established to solve the inhomogeneous
Dirichlet problem on manifolds with boundary and cylindrical ends.

The class of manifolds with boundary and cylindrical ends that we
consider have a product structure at infinity (including the boundary and
the metric). It is possible to relax somewhat these conditions, but for
simplicity we do not address this technical question in this paper.

\begin{definition}\label{def.m.b.c.e}\ Let $N$ be a Riemannian manifold with
boundary $\pa N$. We shall say that $N$ is a {\em manifold with
boundary and cylindrical ends} if there exists an open subset $V$ of
$N$ isometric to $(-\infty, 0) \times X$, where $X$ is a compact
manifold with boundary, such that $N \smallsetminus V$ is compact.
\end{definition}

\begin{lemma}\label{lemma.m.b.c.e}\ Let $N$ be a Riemannian manifold with
boundary $\pa N$. Then $N$ is a manifold with boundary and
cylindrical ends if, and only if, there exists a manifold with
cylindrical ends (without boundary) $M$ with a standard decomposition
$$
    M = M_1 \cup \bigl(\pa M_1 \times (-\infty, 0] \bigr)
$$
and containing $N$ such that
$$
    N \cap \bigl(\pa M_1 \times (-\infty, 0] \bigr) = X \times
    (-\infty,0],
$$
for some compact manifold with boundary $X \subset \pa M_1$.
\end{lemma}

\begin{proof}\ If the metric on $N$ is a product metric on a
tubular neighborhood of $\pa N$, then we can take $M := N \cup (-N)$
to be the double of $N$. The general case can be reduced to this one,
because any metric on $N$ is equivalent to a product metric in a small
tubular neighborhood of $\pa N$.
\end{proof}

Let $M = M_1 \cup (\pa M_1 \times (-\infty, 0])$ be a manifold with
cylindrical ends. The transformation
$$
    (-\infty, -1] \ni x \to t := x^{-1} \in [-1,0)
$$
then extends to a diffeomorphism $\psi$ between $M$ and the interior
$M_0 := M_1 \smallsetminus \pa M_1$ of $M_1$:
\begin{equation}
    \psi : M \to M_0 := M_1 \smallsetminus \pa M_1.
\end{equation}
If $N \subset M$ is a manifold with boundary and cylindrical ends, as
in Lemma \ref{lemma.m.b.c.e}, then the above diffeomorphism will map
$N$ to a subset $N_0 \subset M_0$, whose closure $N_1$ is a {\em
compact} manifold with corners of codimension at most two, $N_1
\subset M_1$. We can identify $N_1$ with the disjoint union $N_0 \cup
X$, if $X$ is as in the definition above.

We shall fix $N \subset M$ as above in what follows.
We define then $H^s(N)$ to be the space of restrictions to the interior of
$N$ of distributions $u \in H^s(M)$. Recall that the main goal of this paper
is to  prove that the map
\begin{equation}\label{eq.nh.Dirichlet}
    H^s(N) \ni u \to ((\Delta_N + V)u , u \vert_{\pa N}) \in
    H^{s-2}(N) \oplus H^{s-1/2}(\pa N)
\end{equation}
is an isomorphism for $s > 1/2$, where $V \ge 0$ a smooth function
that is asymptotically translation invariant in a neighborhood of
infinity (that is $V \in \APS{0}$).

We shall use the results of the previous subsections for the
particular case when $Z = \pa N$.

\begin{proposition}\label{prop.4.2}\ Assume that the potential $V$
is chosen so that $V$ is not identically zero on $\pa M_1\setminus\bar{X}$.
Then the map $-\frac{1}{2}I+K^*:L^2(\pa N)\to L^2(\pa N)$ is injective.
\end{proposition}

\begin{proof}\
Just follow word for word \cite[Proposition 4.1]{MTJFA}.
\end{proof}

Note that our signs are opposite to those in \cite{MTJFA} or
\cite{TaylorC}, because we use the definition that makes the Laplace operator
is {\em positive}.

To prove the Fredholm property of the operators $\frac{1}{2}I+K$ and
$\frac{1}{2}I+K^*$, we need to slightly change the corresponding
argument in \cite{MTJFA}.

\begin{proposition}\label{prop.4.4}\ Retain the same assumptions as in
Proposition~\ref{prop.4.2}. Then the operator

\begin{equation}\label{operator}
    -\ha I+K:L^2(\pa N)\longrightarrow L^2(\pa N)
\end{equation}

\noindent is Fredholm of index zero.
\end{proposition}

\begin{proof}\ The above proposition is known when
$M$ is compact (see \cite[Corollary 4.5]{MTJFA}).

To check that it is Fredholm, we shall rely on $(iv)$ in
Theorem~\ref{theorem.classical2} which, in view of Proposition~\ref{prop.sk},
(\ref{fourier}), and (\ref{multiplier}), amounts to studying the associated
indicial family.

Let $W := V_{\pa M_1}$ and $T=(\Delta + V)^{-1}\pa_\nu^*$. Recall that
$K := T_{\pa N}$ and that $\pa N \sim \pa X\times(-\infty,0]$ in a
neighborhood of infinity. Then Proposition~\ref{prop.comp} gives

\begin{equation}\label{ktau}
\hat{K}(\tau)=\widehat{T_{\pa N}}(\tau)=[\hat T(\tau)]_{\pa X}
=[(\Delta_{\pa M_1}+\tau^2+W)^{-1}\pa^*_{\nu}]_{\pa X}=K_\tau,
\end{equation}

\noindent where $K_\tau$ is the double layer potential operator
associated with the perturbed Laplacian $\Delta_{\pa M_1}+\tau^2+W$ on
$\pa X$ (cf. the discussion in \S{3.3}).
Let $f_{\tau}(x)$ be the Fourier transform in the $t$-variable of $f(x,t)$
($t \in \RR$). In light of this and (\ref{multiplier}), there remains to
prove that the map

\begin{equation}\label{map}
L^2(\pa X\times\RR)\ni f(x,t)\mapsto\mathcal{F}^{-1} \Bigl[\bigl(-\ha
I+K_\tau\bigr)\hat{f}_\tau(x)\Bigr](t) \in L^2(\pa X\times\RR)
\end{equation}

\noindent is an isomorphism. To see this, let $g\in L^2(\pa X\times\RR)$
be arbitrary and, for each $\tau\in\RR$, introduce
$h_\tau:=(-\frac12 I+K_\tau)^{-1}\hat{g}_\tau$.
{}From Proposition~\ref{prop.3.3.1} (utilized for $\Omega:=\pa
M_1\setminus\bar{X}$, which accounts for a change in sign as far as
the coefficient $1/2$ is concerned), it follows that this is
meaningful, $h_\tau\in L^2(\pa X)$ and

\begin{equation}\label{estimateongtau}
    \|h_\tau\|_{L^2(\pa X)}\leq C\|\hat{g}_\tau\|_{L^2(\pa
    X)},\quad \mbox{ uniformly for }\tau\in\RR.
\end{equation}

\noindent If we now set $h(x,t):=\mathcal{F}^{-1}(h_\tau(x))(t)$ then,
thanks to \eqref{estimateongtau} and Plancherel's formula,

\begin{multline}\label{estimateong}
    \int_{\pa X}\int_{\RR}|h(x,t)|^2\,dtd\sigma_x = \int_{\pa
    X}\int_{\RR}|h_\tau(x)|^2\,d\tau d\sigma_x =
    \int_{\RR}\|h_\tau\|^2_{L^2(\pa X)}\,d\tau \\ \leq
    C\int_{\RR}\|\hat{g}_\tau\|^2_{L^2(\pa X)}\,d\tau= C\int_{\pa
    X}\int_{\RR}|g(x,t)|^2\,dtd\sigma_x.
\end{multline}

\noindent That is, $h \in L^2(\pa X\times\RR)$ and $\|h\|_{L^2(\pa
X\times\RR)}\leq C\|g\|_{L^2(\pa X\times\RR)}$.  Furthermore,

\begin{eqnarray} %\label{finalstep}
    \mathcal{F}^{-1}[(-\ha I + K_\tau)\hat{h}_\tau(x)](t) &=&
    \mathcal{F}^{-1}[(-\ha I + K_\tau) h_\tau(x)](t)\nonumber\\ &=&
    \mathcal{F}^{-1}(\hat{g}_\tau)(x)\\ &=& g(x,t)\nonumber
\end{eqnarray}

\noindent which proves that the map (\ref{map}) is onto.  The fact
that \eqref{map} is also one-to-one, follows more or less directly
from the analogue of (\ref{eq.3.3.25}) in our context.

Thus, at this stage, we may conclude that (\ref{operator}) is indeed a
Fredholm operator; there remains to compute its index.  To set the
stage, let us observe that Proposition~\ref{prop.4.2} and duality can
now be used to justify that

\begin{equation}\label{onto}
    -\ha I+K:L^2(\pa N)\longrightarrow L^2(\pa N)\mbox{ is onto}.
\end{equation}

\noindent Next, so we claim,

\begin{equation}\label{onto2}
    -\ha I+K : H^1(\pa N) \longrightarrow H^1(\pa N) \mbox{ is
         Fredholm and onto}
\end{equation}

\noindent as well. Indeed, since $K\in\Psi^{-1}_{\rm ai}(\pa N)$, it
follows that for each $s$,

\begin{equation}\label{smoothing}
    f\in H^s(\pa N)\,\,\&\,\,(-\ha I+K)f\in H^{s+1}(\pa
    N)\Longrightarrow f\in H^{s+1}(\pa N).
\end{equation}

\noindent In concert with (\ref{onto}), this shows that the operator in
(\ref{onto2}) is onto. Also, since

\begin{equation}\label{ker}
    {\rm dim}\,{\rm Ker}\,\Bigl(-\ha I+K;H^1(\pa N)\Bigr)\leq {\rm
    dim}\,{\rm Ker}\,\Bigl(-\ha I+K;L^2(\pa N)\Bigr)<+\infty,
\end{equation}

\noindent the claim (\ref{onto2}) is proved. In particular,

\begin{equation}\label{lesszero}
    {\rm index}\,\Bigl(-\ha I+K;L^2(\pa N)\Bigr)\leq
    0\quad\mbox{and}\quad {\rm index}\,\Bigl(-\ha I+K;H^1(\pa
    N)\Bigr)\leq 0.
\end{equation}

We now take an important step by proving that

\begin{equation}\label{operator2}
    S:L^2(\pa N) \longrightarrow H^1(\pa N) \mbox{ is Fredholm}.
\end{equation}

\noindent (Later on we shall prove that this operator is in fact invertible).
This task is accomplished much as before, i.e. by relying
on Theorem~\ref{theorem.classical2} and Proposition~\ref{prop.3.3.1}, and
we only sketch the main steps. First, as pointed out in
Proposition~\ref{prop.sk}, $S$ is elliptic. Second, the first estimate
in Proposition~\ref{prop.3.3.1} eventually allows us to conclude that the
assignment

\begin{equation*}\label{map2}
    L^2(\pa X\times\RR)\ni f(x,t)\mapsto
    \mathcal{F}^{-1}\Bigl[S_\tau\hat{f}_\tau(x)\Bigr](t)
    \in H^1(\pa X\times\RR)
\end{equation*}

\noindent is an isomorphism, concluding the proof of the claim
(\ref{operator2}).

Having dealt with (\ref{operator2}), we next invoke an
intertwining identity, to the effect that

\begin{equation*}\label{intertwin}
    \bigl(-\ha I+K\bigr)S=S\bigl(-\ha I+K^*\bigr).
\end{equation*}

\noindent This can be seen by starting with Green's formula
$u=\Dl(u|_{\partial N})-\Sl(\partial_\nu u)$ written for the harmonic
function $u:=\Sl(f)$, and then using the jump-relations deduced in
Theorems~\ref{theorem.map.slp}-\ref{theorem.map.dlp}.
The identity (\ref{intertwin}) allows us to obtain

\begin{multline*} %\label{index}
{\rm index}\,\bigl(-\ha I+K;H^1(\pa N)\bigr)
 = {\rm index}\,\bigl(-\ha I+K^*;L^2(\pa N)\bigr)\\
 = -{\rm index}\,\bigl(-\ha I+K;L^2(\pa N)\bigr).
\end{multline*}

\noindent {}From this and (\ref{lesszero}) we may finally conclude
that the operator (\ref{operator}) has index zero, as desired.
\end{proof}

\begin{corollary}\label{cor.iso}\
Let $V$ be as before. Then the operator

\begin{equation*}\label{konhs}
    -\ha I+K:H^s(\pa N)\longrightarrow H^s(\pa N)
\end{equation*}

\noindent is invertible for each $s\in\RR$.
\end{corollary}

\begin{proof}\
To begin with, the case $s=0$ is easily proved by putting together the
above two propositions. In particular, the operator $-\ha I+K:H^s(\pa
N)\longrightarrow H^s(\pa N)$, in the statement of this corollary, is
injective for each $s\geq 0$. Since the fact that this operator is
also surjective is a consequence of the corresponding claim in the
case $s=0$ and the smoothing property (\ref{smoothing}), the desired
conclusion follows for $s\geq 0$.  As for the case $s<0$, a similar
reasoning shows that

\begin{equation}\label{kstaronhs}
-\ha I+K^*:H^{-s}(\pa N)\longrightarrow H^{-s}(\pa N)
\end{equation}

\noindent is invertible for each $s<0$. This and duality then
yield the invertibility of $-\ha I+K:H^s(\pa N)\longrightarrow H^s(\pa N)$
for $s<0$, as wanted.
\end{proof}

Another proof of the above result can be obtained from Theorem
\ref{theorem.sp.i}, for the case $m=0$, the ``easy one.''

Recall that $H^{s}(N)$ is the space of restrictions of distributions in
$H^{s}(M)$ to the interior of $N$. After these preliminaries, we are finally
in a position to discuss the following basic result.

\begin{theorem}\label{theorem.Dirichlet}\
Let $V \in \APS{0}$ be a smooth positive function. For any $s>0$ and
any $f \in H^s(\pa N)$, there exists a unique function $u \in
H^{s+1/2}(N)$ such that $u\vert_{\pa N}=f$ and $(\Delta_N + V)u=0$.
\end{theorem}

\begin{proof}\
Extend first $V$ to a smooth positive function in $\APS{0}$ (that is,
asymptotically translation invariant in a neighborhood of infinity)
which is not identically zero on the complement of $N$.  The
conclusion in Corollary~\ref{cor.iso} will hold for this
function. First we claim that

\begin{equation}\label{Dmapping}
    \Dl:H^s(\pa N)\longrightarrow H^{s+1/2}(N),\qquad s\in\RR,
\end{equation}

\noindent is well-defined and bounded. Indeed, if $s<0$, then this is
a consequence of the implication

\begin{equation}\label{regularity}
    f\in H^s(\pa N),\,\,s<0\,\,\Longrightarrow f\otimes\delta_{\pa
        N}\in H^{s-1/2}(\pa N)
\end{equation}

\noindent along with the factorization $\Dl(g)=(\Delta +
V)^{-1}\partial_\nu^*(g\otimes\delta_{\pa N})$.  For $s=0$, one can
employ the techniques of \cite{MMT}.  The case $s>0$ then follows
inductively from what we have proved so far with the aid of a
commutator identity which essentially reads $\nabla\Dl
f=\Dl(\nabla_{tan}f)+\mbox{lower order terms}$; see (8.19) in
\cite{MT2000} as well as (6.17) in \cite{MMT}.

Having disposed off (\ref{Dmapping}) the existence part in the theorem
is then easily addressed. Specifically, if $s>0$, consider $g:=(-\ha
I+K)^{-1}f\in H^s(\pa N)$ and then set $u:=\Dl(g)\in H^{s+1/2}(N)$ by
(\ref{Dmapping}).

To prove uniqueness, assume that $u\in H^{s+1/2}(N)$ is a null
solution for the Dirichlet problem in $N$. For an arbitrary function
$\varphi\in{\mathcal C}^\infty_c(N)$, let $v$ solve the Dirichlet
problem
$$
    (\Delta_N + V)v = 0, \qquad v \vert_{\partial N} =
    -[(\Delta+V)^{-1}\varphi]|_{\partial N},
$$
and then set $w := v + (\Delta+V)^{-1}\varphi$. It follows that
$(\Delta_N + V) w = \varphi$ in $N$ and $w\vert_{\partial
N}=0$. Consequently, Green's formula gives
$$
    (u, \varphi) = (u, (\Delta_N + V)w)=((\Delta_N + V)u, w)=0
$$
since $u\vert_{\partial N} = w\vert_{\partial N}=0$. Since $\varphi$
is arbitrary, this forces $u=0$ in $N$ as desired.
\end{proof}

We are now ready to prove Theorem~\ref{theorem.main} which, for
the convenience of the reader, we restate below.

\begin{theorem}\label{theorem.main.again}\
Let $N$ be a manifold with boundary and cylindrical ends and $V \ge 0$
be a smooth functions that is asymptotically translation invariant in
a neighborhood of infinity. Then
$$
    H^{s}(N) \ni u \to \tilde \Delta_N (u) := ((\Delta_N + V)u,
    u\vert_{\pa N}) \in H^{s-2}(N) \oplus H^{s-1/2}(\pa N)
$$
is a continuous bijection, for any $s > 1/2$.
\end{theorem}

\begin{proof}\ First we extend $V$ to $M$, making sure that it is
still $\ge 0$, smooth, and asymptotically translation invariant.  The
continuity of the map $\tilde \Delta_N$ follows from the continuity of
$\Delta_N + V : H^s(N) \to H^{s-2}(N)$ and from the continuity of the
trace map $H^s(N) \to H^{s - 1/2}(\pa N)$.

As before, we fix a potential $V$ which vanishes in a neighborhood of
$N$.  Let $g \in H^{s-2}(N)$ be arbitrary. First extend $g$ to a
distribution (denoted also $g$) in $H^{s-2}(M)$, then set $u_1 =
(\Delta+V)^{-1}g\in H^{s}(M)$ and $f_1 = u_1\vert_{\pa N}\in
H^{s-1/2}$.  Finally, choose $u_2 \in H^{s}(N)$ such that $(\Delta_N +
V)u_2=0$ and $u_2\vert_{\partial N}=f-f_1$. Then $u := u_1 + u_2$
satisfies $(\Delta_N + V)u = g$ and $u \vert_{\pa N} = f$. This proves
the surjectivity of $\tilde\Delta_N$. The injectivity of this map then
follows from the uniqueness part in Theorem~\ref{theorem.Dirichlet}.
\end{proof}

It is likely that some versions of the above two theorems extend
to weighted Sobolev spaces. This will likely requires techniques
similar to those used in \cite{GilMendoza}. In \cite{ScSc},
Schrohe and Schulze have generalized the Boutet de Monvel calculus
to manifolds with boundary and cylindrical ends. With some
additional work, their results can probably be used to prove our
Theorem~\ref{theorem.main} above. Our approach, however, is
shorter and also leads to a characterization of the
Dirichlet-to-Neumann boundary map, Theorem \ref{theorem.DtoN}. It
is worth pointing out that our methods can also handle non-smooth
structures (cf. \S{4}) and seem amenable to other basic problems
of mathematical physics in non-compact manifolds (such as
Maxwell's equations in infinite cylinders). We hope to return to
these issues at a later time.

\subsection{The Dirichlet-to-Neumann map}
Theorem~\ref{theorem.Dirichlet} allows us to define the
Dirichlet-to-Neumann map $\mathcal N$
$$
    \mathcal N(f) = (\pa_\nu u)_+
$$
for $f \in L^2(\pa N)$ and $u$ solution of $(\Delta_N + V)u =0$, 
$u_+:= u\vert_{\pa N} = f$.

\begin{theorem}\label{theorem.DtoN}\ Let $N$ be a manifold with boundary
and cylindrical ends. Then the operator $S : H^s(\pa N) \to
H^{s+1}(\pa N)$ of Equation~\eqref{eq.def.S} is invertible for any $s$
and $(\ha I + K^*)S^{-1} = \mathcal N$, the ``Dirichlet-to-Neumann
map.'' In particular, $\mathcal N \in \APSN{1}$.
\end{theorem}

\begin{proof}\ The operator $S$ is elliptic by Proposition \ref{prop.sk}.
For further reference, let us note here that

\begin{equation}\label{S-smoothing}
    f\in H^s(\pa N) \, \, \& \, \, Sf\in H^{s+1}(\pa N)
    \Longrightarrow f\in H^{s+1}(\pa N),
\end{equation}

\noindent by elliptic regularity.

Next, using the notation of Proposition~\ref{prop.3.3.1}, we have
$\hat S(\tau) = S_\tau$. By the results of the same proposition, $\hat
S(\tau)$ is invertible for any $\tau$, and the norm of the inverse is
uniformly bounded (this can be proved also by using the results of
\cite{MTJFA} or \cite{MMT} and the estimates in \cite{Shubin}).
Consequently, $S : H^s(\pa N) \to H^{s + 1}(\pa N)$ is Fredholm
(cf. Theorem~\ref{theorem.classical2}).

Checking that $S$ is injective when $s=0$ is done much as in the last
part of \S{6} in \cite{MTJFA}. In short, the idea is as follows.
Assume that $Sf = 0$ for some $f \in H^s(\pa N)$ and let $N \subset
M$, where $M$ is a manifold with cylindrical ends without boundary, as
in Lemma~\ref{lemma.m.b.c.e}. Then $u:=\Sl(f)$ satisfies $(\Delta +
V)u = 0$ on $M \smallsetminus \pa N$, and
$$
    u_{\pa N} = u_+ = u_- = Sf = 0.
$$
Furthermore, thanks to (\ref{S-smoothing}), (\ref{regularity}) and the
factorization $\Sl(f) = (\Delta + V)^{-1}(f \otimes \delta_{\partial
N})$, the function $u$ is sufficiently regular so that (the uniqueness
part in) Theorem~\ref{theorem.Dirichlet} holds both in $N$ and in
$M\setminus N$.  Hence, by Theorem~\ref{theorem.map.slp},
$$
    f = (\pa_\nu u)_+ -(\pa_\nu u)_- = 0,
$$
as desired. Thus, $S : H^s(\pa N) \to H^{s + 1}(\pa N)$ is injective,
first for $s\geq 0$ (via a simple embedding), then for $s\in\RR$ via
(\ref{S-smoothing}).

Since $S$ is formally self-adjoint, we get that $S$ has also dense
range.  Using now the fact that $S$ is Fredholm, we obtain that $S$ is
bijective, as desired.
\end{proof}

\begin{corollary}\label{cor.DtoN}\
The Cauchy data space
$$
    \{(u\vert_{\pa N}, \pa_\nu u\vert_{\pa N});\, u \in H^s(N)\,,
    \; (\Delta_{N} + V) u = 0 \}
$$
is a closed subspace of $H^{s-1/2}(\pa N) \oplus H^{s-3/2}(\pa N)$ for
any $s > 1/2$.
\end{corollary}

\begin{proof}\ By Theorem \ref{theorem.main.again}, the Cauchy data space
$$
    \mathcal C := \{(u\vert_{\pa N}, \pa_\nu u\vert_{\pa N});\, u
    \in H^s(N)\,,\; (\Delta_{N} + V) u = 0 \}
$$
is given by the graph of $\mathcal N$, namely
$$
    \mathcal C = \Gamma(\mathcal N) : = \{ (f, \mathcal N f), f
    \in H^{s-1/2} \} \subset H^{s-1/2}(\pa N) \oplus H^{s-3/2}(\pa
    N).
$$
Theorem \ref{theorem.DtoN} shows that this space is closed, since
$\mathcal N \in \APSN{1}$ and hence it defines a continuous
(everywhere defined) map $H^{s-1/2}(\pa N) \to H^{s-3/2}(\pa N)$.
\end{proof}

We conclude this section with yet another integral representation
formula for the Dirichlet problem.

\begin{corollary}\label{cor.SSinv}\
Retain the usual set of assumptions. Then, for each $s>0$, the
solution to the boundary problem
$$
    u \in H^{s+1/2}(N), \quad (\Delta_N + V)u=0, \quad u\vert_{\pa
    N} = f\in H^s(\pa N),
$$
(first treated in Theorem~\ref{theorem.Dirichlet}) can also be
expressed in the form
$$
    u=\Sl(S^{-1}f)\mbox{ in }N.
$$
\end{corollary}

\begin{proof}\ The starting point is the claim
(which can be justified in a manner similar to (\ref{Dmapping})) that

\begin{equation}\label{Smapping}
\Sl:H^s(\pa N)\longrightarrow H^{s+3/2}(M),\qquad s\in\RR,
\end{equation}

\noindent is well-defined and bounded. In concert with the fact that
$S : H^s(\pa N) \to H^{s + 1}(\pa N)$ is invertible, this finishes the
proof of the corollary.
\end{proof}

\end{document}